\documentclass[12pt]{amsart}
\usepackage{amsmath,amssymb,mathrsfs,hyperref,color}

\usepackage[nameinlink]{cleveref}
\hypersetup{colorlinks={true},linkcolor={blue},citecolor=green}
\crefformat{equation}{(#2#1#3)}

\usepackage{mathtools}
\usepackage[x11names]{xcolor}
\newtagform{blue}{\color{blue}(}{)}

\makeatletter
\def\@fnsymbol#1{\ensuremath{\ifcase#1\or *\or 1\or 2\or
   3\or 4\or 5\or 6\or 7\or 8\else\@ctrerr\fi}}
\makeatother

\usepackage{paralist}

\usepackage{subfigure}
\usepackage{float}
\usepackage{times}
\usepackage{tikz}
\usepackage{cases}
\usetikzlibrary{intersections}
\usepackage{xcolor}
\usepackage[latin1]{inputenc}
\usepackage{bbm}
\usetikzlibrary{shadings}
\usetikzlibrary[intersections]
\usetikzlibrary{arrows,shapes}

\makeatletter
\def\setliststart#1{\setcounter{\@listctr}{#1}%
  \addtocounter{\@listctr}{-1}}
\makeatother

\makeatletter
\@addtoreset{figure}{section}
\makeatother

\usepackage{calc}
\newtheorem{theorem}{Theorem}[section]
\newtheorem{lemma}[theorem]{Lemma}
\newtheorem{proposition}[theorem]{Proposition}
\newtheorem{corollary}[theorem]{Corollary}
\newtheorem{remarks}[theorem]{Remark}

\numberwithin{equation}{section}

\newcommand{\R}{\mathbb{R}}

\newcommand{\N}{\mathbb{N}}

\newcommand{\PP}{\mathcal{P}}

\DeclareMathOperator*{\supp}{spt}

\DeclareMathOperator*{\esssup}{ess\ sup}
\DeclareMathOperator*{\ddiv}{div}

\DeclareMathOperator*{\eps}{\varepsilon}

\makeatletter
\def\moverlay{\mathpalette\mov@rlay}
\def\mov@rlay#1#2{\leavevmode\vtop{%
   \baselineskip\z@skip \lineskiplimit-\maxdimen
   \ialign{\hfil$\m@th#1##$\hfil\cr#2\crcr}}}
\newcommand{\charfusion}[3][\mathord]{
    #1{\ifx#1\mathop\vphantom{#2}\fi
        \mathpalette\mov@rlay{#2\cr#3}
      }
    \ifx#1\mathop\expandafter\displaylimits\fi}
\makeatother

\title[A singular perturbation problem for mean field games of acceleration]{A singular perturbation problem for mean field games of acceleration: application to mean field games of control}
\author{Cristian Mendico}
\address{Dipartimento di Matematica, Universit\`a degli studi di Roma Tor Vergata -- Via della Ricerca Scientifica 1, 00133 Roma}
\email{mendico@mat.uniroma2.it}

\date{\today}
\thanks{{\it Acknowledgement and funding:} Cristian Mendico was partly supported by Istituto Nazionale di Alta Matematica (GNAMPA 2020 Research Projects) and by the MIUR Excellence Department Project awarded to the Department of Mathematics, University of Roma Tor Vergata, CUP E83C23000330006. The author would like to thanks Pierre Cardaliaguet for his fruitful comments and his careful reading of the manuscript.}

\begin{document}
\usetagform{blue}
\begin{abstract}
We study the singular perturbation problem for mean field game systems with control of acceleration. For such a problem we analyze the behavior of solutions as the acceleration costs vanishes. In this setting the Hamiltonian fails to be strictly convex and coercive  w.r.t. the momentum variable and this creates new issues in the analysis of the problem. We show that the limit system is of MFG type: we first study the convergence to the classical MFG system and, then, by a finer analysis of the Euler-Lagrange flow associated with the control of acceleration we show the convergence to a class of, so-called, MFG of control problems.
\end{abstract}

\maketitle
\tableofcontents
\section{Introduction}

The study of singular perturbation problem of control systems has a long history going back to \cite{bib:ZA, bib:ARG, bib:ZAV} and references therein. Such a problem concerns the analysis of systems where some state variables evolve at a much faster time scale than the others. Generally, the solution of a typical singular perturbation problem leads to the elimination of the fast state variable and, consequently, to a reduction of the dimension of the system. Clearly, the limit problem keeps some informations on the fast part.

Besides classical control systems, other type of singular perturbation  problems have  been studied and we refer, for instance, to homogenization (e.g. \cite{bib:LPV}, \cite{bib:IH}) and the long time behavior (e.g. \cite{bib:FM, bib:FA}, \cite{bib:BJMR}). More recently, such analysis have been extended to the case of differential games (e.g. \cite{bib:ABG, bib:ABE}, \cite{bib:MG}, \cite{bib:KDV}) and of mean field games (MFG) (e.g. \cite{bib:CAR}, \cite{bib:CCMW, bib:CCMW1}, \cite{bib:PC1}, \cite{bib:CDM, bib:SLL}).  Based on this recent literature, in this paper we make a step further. Indeed, the goal of this work is twofold: first, we show a connection between the classical MFG system, where the underlying payoff is a calculus of variation problem, and the MFG with control of acceleration; secondly, we show how a MFG of control system can be recovered from a MFG system with control of acceleration. We will extend such analysis to the study of singular perturbation problems associated with sub-Riemannian structure and MFG defined on such structures in a future work.

We recall that MFG were introduced in  \cite{bib:LL1, bib:LL2, bib:LL3} and \cite{bib:HCM1, bib:HCM2} in order to describe the behavior of Nash equilibria in problems with infinitely many rational agents (we refer to \cite{bib:NC} and references therein for more details). Since these pioneering works the MFG theory has grown very fast:  we refer for instance to the survey papers and the monographs \cite{bib:DEV, bib:BFY, bib:CD1}. The classical MFG system introduced in  \cite{bib:LL1, bib:LL2, bib:LL3} describes systems in which each the typical payoff is represented by deterministic calculus of variation problem. MFG systems with control of acceleration, first introduced in \cite{bib:CM, bib:YA}, describe models where agents control their acceleration and the cost functional to minimize depends on higher order derivatives of admissible trajectories. Such problems naturally appear in the study of agent-based models which describe the collective behavior of various animal populations (e.g. \cite{bib:AC, bib:TBL}) or crowd dynamics (e.g. \cite{bib:CPT, bib:EBA}). In this framework the study of the singular perturbation problem we perform in this paper finds a lot of applications: for instance, such an analysis can be applied to a MFG system of Cucker-Smale type, see for instance \cite{Bardi_2021}, to describe the behavior of a flock in which the control is increasingly cheap.

We describe now the problems we are going to solve in this paper.
\begin{enumerate}
\item {\bf Convergence to the classical MFG system.} 
We study the limit of the solution to the system 
\begin{align}\label{eq:AppepsMFG}
 \begin{cases}
 	-\partial_{t} u^{\eps} +\frac{1}{2\eps}|D_{v}u^{\eps}|^{2} - \langle D_{x}u^{\eps}, v \rangle -L_{0}(x, v, m^{\eps}_{t})= 0, & (t,x,v) \in [0,T] \times \R^{2d}
 	\\
 	\partial_{t}\mu^{\eps}_{t} - \langle D_{x}\mu^{\eps}_{t},v \rangle - \frac{1}{\eps}\ddiv_{v}\left(\mu^{\eps}_{t}D_{v}u^{\eps} \right)=0, &  (t,x,v) \in [0,T] \times \R^{2d}
\\
\mu^{\eps}_{0}=\mu_{0}, \quad u^{\eps}(T,x,v)=g(x,m^{\eps}_{T}), & (x,v) \in \R^{2d}
 \end{cases}	
 \end{align}
as the parameter $\eps$ goes to zero. Heuristically, the state equation associated with the above PDEs system is given by 
\begin{equation}\label{eq:dynamics}
\begin{cases}
\dot x(t)=v(t)
\\
\dot v(t)= \frac{1}{\eps} \alpha(t)
\end{cases}
\end{equation}
where $\alpha : [0,T] \to \R^{d}$ is a measurable control function and, from \cite{bib:CM, bib:YA}, we have that for any $\eps > 0$ a typical player aims to minimize a cost functional of the form
\begin{equation*}
\int_{t}^{T}{\left(\frac{\eps}{2}|\ddot\gamma(s)|^{2} + L_{0}(\gamma(s), \dot\gamma(s), m^{\eps}_{s}) \right)\ ds} + g(\gamma(T), m^{\eps}_{T}). 
\end{equation*}
Moreover, still from \cite{bib:CM, bib:YA},  under suitable assumptions (listed below) on the function $L_{0}$, we have that for any $\eps > 0$ there exists a unique solution $(u^{\eps}, \mu^{\eps})$ to \eqref{eq:AppepsMFG}.

Following the previous considerations on a typical singular perturbation problem,  in case of control of acceleration we expect that the fast variable, in this case the velocity of each player, is eliminated in the limit and all the informations are captured by the behavior of the space variable. Moreover, since the aim of such analysis is to establish a rigorous mathematical connection between the classical MFG system and the MFG system with control of acceleration, system \eqref{eq:AppepsMFG} has such a particular form. Indeed, we observe that the function $L_{0}$ and the terminal costs $g$ only depend on the space marginal of the measures $\{\mu^{\eps}_{t}\}_{t \in [0,T]}$. Such a marginal flow of measures captures the behavior of the fast variable in the limit and it is also the object of investigation in classical MFG since it represents the distribution of players in space at each time $t \in [0,T]$.

\item {\bf Convergence to MFG of control system.} 
 In the second part, we analyze the limit of the solution to the system 
\begin{align}\label{eq:AppepsMFG}
 \begin{cases}
 	-\partial_{t} u^{\eps} +\frac{1}{2\eps}|D_{v}u^{\eps}|^{2} - \langle D_{x}u^{\eps}, v \rangle -L_{0}(x, v, \mu^{\eps}_{t})= 0, & (t,x,v) \in [0,T] \times \R^{2d}
 	\\
 	\partial_{t}\mu^{\eps}_{t} - \langle D_{x}\mu^{\eps}_{t},v \rangle - \frac{1}{\eps}\ddiv_{v}\left(\mu^{\eps}_{t}D_{v}u^{\eps} \right)=0, &  (t,x,v) \in [0,T] \times \R^{2d}
\\
\mu^{\eps}_{0}=\mu_{0}, \quad u^{\eps}(T,x,v)=g(x,m^{\eps}_{T}), & (x,v) \in \R^{2d}
 \end{cases}	
 \end{align}
still as the parameter $\eps$ goes to zero. The main issue here is that both the data $L$ and $g$ depend on $\mu^{\eps}$ and we have to deal with the convergence of the whole measure. Note that, even though the limit control problem does not depend on velocity as a state variable we have that the second marginal of the limit measure, and so the Lagrangian function, still depends on it. For this reason, we expect the limit system to be of mean field game of control type. 
\end{enumerate}

Next, we briefly explain the main result of this work and the method of proof. 
\begin{enumerate}
\item {\bf Towards the classical MFG system.} We prove that $(u^{\eps}, m^{\eps})$, where $m^{\eps}_{t}$ is the space marginal of the solution $\mu^{\eps}_{t}$ for any $t \in [0,T]$, converges (up to subsequence) to a solution $(u^{0}, m^{0})$ of the classical MFG system 
\begin{align}\label{eq:introappMFG}
	\begin{cases}
		(i)\,\, -\partial_{t} u^{0}(t,x) + H_{0}(x, D_{x}u^{0}(t,x), m^{0}_{t})=0, & \quad (t,x) \in [0,T] \times \R^{d}
		\\
		(ii)\,\, \partial_{t}  m^{0}_{t} - \ddiv\Big( m^{0}_{t}D_{p}H_{0}(x, D_{x}u^{0}(t,x), m^{0}_{t}) \Big)=0, & \quad (t,x) \in [0,T] \times \R^{d}
		\\
		m^{0}_{0}= \pi_{1} \sharp \mu_{0},\,\, u^{0}(T,x)=g(x,m^{0}_{T}), & \quad x \in \R^{d}
	\end{cases}
	\end{align}
	where $H_{0}:\R^d \times \R^d \to \R$ is the Legendre Transform of the function $L_0$. Observe that, the main difference between our result and the existing one concerning the homogenization problem in MFG (\cite{bib:CDM}, \cite{bib:SLL}, \cite{Bardi_2021}) is that the limit system is still of MFG type. Indeed, in \cite{bib:CDM}, \cite{bib:SLL} and \cite{Bardi_2021} it has been proved that in the limit the MFG structure of the problem is lost and, in particular, an explicit example of MFG system with potential coupling function is constructed in \cite{bib:SLL}. 
	
	In order to prove our first main convergence result, we begin by showing that $u^{\eps}$ is equibounded and $m^{\eps}$ is tight (see \Cref{lem:Mequibound} and \Cref{thm:spacetightness}). Thus, as a first consequence we get that, up to a subsequence, there exists $m^{0} \in C([0,T]; \PP_{1}(\R^{d})$ such that $m^{\eps} \to m^{0}$ in $C([0,T]; \PP_{1}(\R^{d}))$. Then, we proceed with the analysis of the value function $u^{\eps}$: we show that $u^{\eps}(t,\cdot,v)$ is equi-Lipschitz continuous, $u^{\eps}(\cdot,x,v)$ is equicontinuous and $u^{\eps}(t,x,\cdot)$ has decreasing oscillation w.r.t. $\eps$ (see \Cref{lem:oscillationvel} and \Cref{prop:Mequicontinuity}). We finally address the locally uniform convergence of $u^{\eps}$, showing that there exists a subsequence $\eps_{k} \downarrow 0$ such that $(u^{\eps_{k}}, m^{\eps_{k}})$ converges to a solution $(u^{0}, m^{0})$ of \eqref{eq:introappMFG} (see \Cref{thm:Meps}, \Cref{prop:convminim} and \Cref{cor:limiteq}). The main issues in proving the above results are due to the lack of strict convexity and the lack of coercivity of the Hamiltonian in system \eqref{eq:AppepsMFG}. The technic we use to study our singular perturbation problem is a combination of variational approach to Hamilton-Jacobi equation and optimal transport in order to overcome the issues mentioned above. 
\item {\bf Towards MFG of control system.} Just for simplicity of notation, we restrict the attention to a Lagrangian of the form 
	\begin{equation*}
	L(x,v,w, m)= \frac{1}{2}|w|^2 +\frac{1}{2}|v|^2 + L_{0} (x, m). 
	\end{equation*} 
In this setting, we prove that $(u^{\eps}, \mu^{\eps})$ converges (up to subsequence) to a solution $(u^{0}, \mu^{0})$ to the MFG of control system 
\begin{align}\label{eq:introMFG}
	\begin{cases}
		(i)\,\, -\partial_{t} u^{0}(t,x) + \frac{1}{2} |D_{x}u^{0}(t,x)|^{2} - L_{0}(x,\mu^{0}_{t})=0, & \quad (t,x) \in [0,T] \times \R^{d}
		\\
		(ii)\,\, \partial_{t}  m^{0}_{t} - \ddiv\big( m^{0}_{t} D_{x}u^{0}(t,x) \big)=0, & \quad (t,x) \in [0,T] \times \R^{d}
		\\
		(iii)\,\, \mu^{0}_{t} = (\text{Id}(\cdot), Du^{0}(t, \cdot)) \sharp m^{0}_{t}
		\\
		\mu^{0}_{0}= \mu_{0},\,\, u^{0}(T,x)=g(x, \mu^{0}_{T}), & \quad x \in \R^{d}
	\end{cases}
	\end{align}
	where $H_{0}:\R^d \times \R^d \to \R$ is the Legendre Transform of the function $L_0$, $m^{0}_{t} = \pi_1 \sharp \mu^{0}_{t}$ and $\text{Id}(\cdot)$ denotes the identity function. As observed before, the main difference with the previous study is the convergence of the whole measure $\mu^{\eps}$ which requires a finer study of the Euler-Lagrange flow associated with the problem of control of acceleration. We observe that equations $(i)$, $(ii)$ are in common with system \eqref{eq:introappMFG} and they differ only in the measure argument of the function $L_{0}$. However, system \eqref{eq:introMFG} has a third equation, $(iii)$, which describes the evolution of the flow $\{\mu^{0}_{t}\}_{t \in [0,T]}$: the second marginal, that is the one w.r.t. the velocity variable, is given by the push-forward of the optimal feedback function $Du^{0}$ by the first marginal $\{m^{0}_{t}\}_{t \in [0,T]}$. Heuristically, such an equation describes the evolution of the density distribution of controls w.r.t. the state of a typical player. For this reason system \eqref{eq:introMFG} is called MFG system of control.  
	
In conclusion, we stress that the result can be generalized to any Lagrangian following the same arguments but with heavy notation that leads to an hard presentation of the ideas. 
\end{enumerate}

\medskip
\noindent The paper is organized as follows. 
In \Cref{sec:preliminaries}  we fix the notation that will be used throughout the paper and we recall the main definitions and results from measure theory. In \Cref{sec:setting} we introduce the MFG system associated with the singular perturbation problem, we give the standing assumptions on the data and finally we state the main results (\Cref{thm:main1} and \Cref{thm:main2}). \Cref{sec:mainresult} and \Cref{sec:proof2} are devoted to the proofs of preliminary results needed to demonstrate \Cref{thm:main1} and \Cref{thm:main2}, respectively.

 \medskip\medskip
 
 \noindent{\bf Statements and Declarations:} There are no associated data. \\ 
 \noindent{\bf Statements and Declarations:} There are no conflict of interests.

\section{Notations and preliminaries}
\label{sec:preliminaries}
\subsection{Notation}
We write below a list of symbols used throughout this paper.
\begin{itemize}
	\item Denote by $\mathbb{N}$ the set of positive integers, by $\mathbb{R}^d$ the $d$-dimensional real Euclidean space,  by $\langle\cdot,\cdot\rangle$ the Euclidean scalar product, by $|\cdot|$ the usual norm in $\mathbb{R}^d$, and by $B_{R}$ the open ball with center $0$ and radius $R$.

\item For a Lebesgue-measurable subset $A$ of $\mathbb{R}^d$, we let $\mathcal{L}^{d}(A)$ be the $d$-dimensional Lebesgue measure of $A$ and  $\mathbf{1}_{A}:\mathbb{R}^n\rightarrow \{0,1\}$ be the characteristic function of $A$, i.e.,
\begin{align*}
\mathbf{1}_{A}(x)=
\begin{cases}
1  \ \ \ &x\in A,\\
0 &x \not\in A.
\end{cases}
\end{align*} 
We denote by $L^p(A)$ (for $1\leq p\leq \infty$) the space of Lebesgue-measurable functions $f$ with $\|f\|_{p,A}<\infty$, where   
\begin{align*}
& \|f\|_{\infty, A}:=\esssup_{x \in A} |f(x)|,
\\& \|f\|_{p,A}:=\left(\int_{A}|f|^{p}\ dx\right)^{\frac{1}{p}}, \quad 1\leq p<\infty.
\end{align*}
For brevity, $\|f\|_{\infty}$  and $\|f\|_{p}$ stand for  $\|f\|_{\infty,\mathbb{R}^d}$ and  $\|f\|_{p,\mathbb{R}^d}$ respectively.



%

\item $C_b(\mathbb{R}^d)$ stands for the function space of bounded uniformly  continuous functions on $\mathbb{R}^d$. $C^{2}_{b}(\mathbb{R}^{d})$ stands for the space of bounded functions on $\mathbb{R}^d$ with bounded uniformly continuous first and second derivatives. 
$C^k(\mathbb{R}^{d})$ ($k\in\mathbb{N}$) stands for the function space of $k$-times continuously differentiable functions on $\mathbb{R}^d$, and $C^\infty(\mathbb{R}^{d}):=\cap_{k=0}^\infty C^k(\mathbb{R}^{d})$. 
 $C_c^\infty(\mathbb{R}^{d})$ stands for the space of functions in $C^\infty(\mathbb{R}^{d})$ with compact support. Let $a<b\in\mathbb{R}$.
  $AC([a,b];\mathbb{R}^d)$ denotes the space of absolutely continuous maps $[a,b]\to \mathbb{R}^d$.
  
  \item For $f \in C^{1}(\mathbb{R}^{d})$, the gradient of $f$ is denoted by $Df=(D_{x_{1}}f, ..., D_{x_{n}}f)$, where $D_{x_{i}}f=\frac{\partial f}{\partial x_{i}}$, $i=1,2,\cdots,d$.
Let $k$ be a nonnegative integer and let $\alpha=(\alpha_1,\cdots,\alpha_d)$ be a multiindex of order $k$, i.e., $k=|\alpha|=\alpha_1+\cdots +\alpha_d$ , where each component $\alpha_i$ is a nonnegative integer.   For $f \in C^{k}(\mathbb{R}^{d})$,
define $D^{\alpha}f:= D_{x_{1}}^{\alpha_{1}} \cdot\cdot\cdot D^{\alpha_{d}}_{x_{d}}f$. 


\end{itemize}

\subsection{The Wasserstein spaces}
We recall here the notations and definitions of Wasserstein spaces and Wasserstein distance, for more details we refer to \cite{bib:CV, bib:AGS}. 

Let $(X,{\bf d})$ be a metric space (in the paper, we use $X= \R^d$ or $X= \R^d\times \R^d$). 
Denote by $\mathcal{B}(X)$ the  Borel $\sigma$-algebra on $X$ and by $\mathcal{P}(X)$ the space of Borel probability measures on $X$.
The support of a measure $\mu \in \mathcal{P}(X)$, denoted by $\supp(\mu)$, is the closed set defined by
\begin{equation*}
\supp (\mu) := \Big \{x \in X: \mu(V_x)>0\ \text{for each open neighborhood $V_x$ of $x$}\Big\}.
\end{equation*}
We say that a sequence $\{\mu_k\}_{k\in\mathbb{N}}\subset \mathcal{P}(X)$ is weakly-$*$ convergent to $\mu \in \mathcal{P}(X)$, denoted by
$\mu_k \stackrel{w^*}{\longrightarrow}\mu$,
  if
\begin{equation*}
\lim_{n\rightarrow \infty} \int_{X} f(x)\,d\mu_n(x)=\int_{X} f(x) \,d\mu(x), \quad  \forall f \in C_b(X).
\end{equation*}

For $p\in[1,+\infty)$, the Wasserstein space of order $p$ is defined as
\begin{equation*}
\mathcal{P}_p(X):=\left\{m\in\mathcal{P}(X): \int_{X} d(x_0,x)^p\,dm(x) <+\infty\right\},
\end{equation*}
for some (and thus all) $x_0 \in X$.
Given any two measures $m$ and $m^{\prime}$ in $\mathcal{P}_p(X)$,  define
\begin{equation}\label{def.transportplan}
\Pi(m,m'):=\Big\{\lambda\in\mathcal{P}(X \times X): \lambda(A\times X)=m(A),\ \lambda(X \times A)=m'(A),\ \forall A\in \mathcal{B}(X)\Big\}.
\end{equation}
The Wasserstein distance of order $p$ between $m$ and $m'$ is defined by
    \begin{equation*}\label{dis1}
          d_p(m,m')=\inf_{\lambda \in\Pi(m,m')}\left(\int_{X\times X}d(x,y)^p\,d\lambda(x,y) \right)^{1/p}.
    \end{equation*}
    The distance $d_1$ is also commonly called the Kantorovich-Rubinstein distance and can be characterized by a useful duality formula (see, for instance, \cite{bib:CV})  as follows
\begin{equation}\label{eq:100}
d_1(m,m')=\sup\left\{\int_{X} f(x)\,dm(x)-\int_{X} f(x)\,dm'(x) \ |\ f:X\rightarrow\mathbb{R} \ \ \text{is}\ 1\text{-Lipschitz}\right\},
\end{equation}
for all $m$, $m'\in\mathcal{P}_1(X)$.

Let $X_{1}$, $X_{2}$ be metric spaces, let $\mu \in \PP(X_{1})$ and let $f: X_{1} \to X_{2}$ be a $\mu$ measurable map. Then, we denote by $f \sharp \mu \in \PP(X_{2})$ the push-forward of $\mu$ through $f$ defined by
\begin{equation*}
f \sharp \mu(B):= \mu(f^{-1}(B)), \quad \forall\ B \in \mathcal{B}(X_{2}).	
\end{equation*}
More generally, in integral form, it reads as
\begin{equation*}
\int_{X_{1}}{\varphi(f(x))\ \mu(dx)} = \int_{X_{2}}{\varphi(y)\ f \sharp \mu(dy)}.	
\end{equation*}

	\section{Setting and main results}
	\label{sec:setting}

	\subsection{Convergence to classical MFG system}
	Let $L_{0}: \R^{2d} \times \PP_{1}(\R^{d}) \to \R$ satisfy the following. 
	 	\begin{itemize}
		\item[{\bf (M1)}] $L_{0}$ is continuous w.r.t. all variables and for any $m \in \PP_{1}(\R^{d})$ the map $(x,v) \mapsto L_{0}(x,v,m)$ belongs to $C^{1}(\R^{d})$.
		\item[{\bf (M2)}] There exists $M_{0} > 0$ such that for any $(x,v,m) \in \R^{2d} \times \PP_{1}(\R^{d})$
		\begin{align}
		D^{2}_{v}L_{0}(x,v,m) \geq\ & \frac{1}{M_{0}}\text{Id},\label{eq:convexity}
			\\
        	|D_{x}L_{0}(x,v,m)| \leq\ & M_{0}\big(1+|v|^{2} \big), &\label{eq:MLspacegradient}
        	\\
        	|D_{v}L_{0}(x,v,m)| \leq\ & M_{0}\big( 1 + |v|\big). \label{eq:MLvelgradient}
        \end{align}
       \item[{\bf (M3)}] There exist two moduli $\theta: \R_{+} \to \R_{+}$ and $\omega_{0}: \R_{+} \to \R_{+}$ such that
        \begin{equation*}
        |L_{0}(x,v,m_{1})-L_{0}(x,v,m_{2})| \leq \theta(|x|) \omega_{0}(d_{1}(m_{1},m_{2})),	
        \end{equation*}
for any $(x,v) \in \R^{2d}$ and $m_{1}$, $m_{2} \in \PP_{1}(\R^{d})$.
	\end{itemize}
	
	Observe that from {\bf (M2)} one easily obtain
	\begin{equation} \label{eq:Mcontrolgrowth}
	\frac{1}{M_{0}}|v|^{2}-M_{0} \leq L_{0}(x,v,m) \leq\  M_{0}(1+|v|^{2}),
	\end{equation}
	  and, without loss of generality,  $L_{0}(x,v,m) \geq 0$ for any $(x,v,m) \in \R^{2d} \times \PP_{1}(\R^{d})$. Let $H_{0}$ be the Legendre Transform of the function $L_{0}$, i.e., 
\begin{equation*}
H_{0}(x,p,m)=\sup_{ v \in \R^{d}} \big\{-\langle p,v \rangle - L_{0}(x,v,m) \big\}.	
\end{equation*}

We consider the MFG system  
	\begin{align}\label{eq:epsMFG}
 \begin{cases}
 	-\partial_{t} u^{\eps} +\frac{1}{2\eps}|D_{v}u^{\eps}|^{2} - \langle D_{x}u^{\eps}, v \rangle -L_{0}(x, v, m^{\eps}_{t})= 0, & (t,x,v) \in [0,T] \times \R^{2d}
 	\\
 	\partial_{t}\mu^{\eps}_{t} - \langle D_{x}\mu^{\eps}_{t},v \rangle - \frac{1}{\eps}\ddiv_{v}\left(\mu^{\eps}_{t}D_{v}u^{\eps} \right)=0, &  (t,x,v) \in [0,T] \times \R^{2d}
\\
\mu^{\eps}_{0}=\mu_{0}, \quad u^{\eps}(T,x,v)=g(x,m^{\eps}_{T}), & (x,v) \in \R^{2d}
 \end{cases}	
 \end{align}
 where  $m^{\eps}_{t}=\pi_{1} \sharp \mu^{\eps}_{t}$ and $\pi_{1}: \R^{2d} \to \R^{d}$ denotes the projection onto the first factor, i.e., $\pi_{1} (x,v)= x$. We assume the following on the boundary data of the system. 
 \begin{itemize}
 \item[{\bf (BC1)}] The measure $\mu_{0} \in \PP(\R^{2d})$ is absolutely continuous w.r.t. Lebesgue measure, we still denote by $\mu_{0}$ its density, and it has compact support. 
 \item[{\bf (BC2)}] The terminal costs $g(\cdot, m)$ belongs to $C^{1}_{b}(\R^{d})$ such that $M_{0} \geq \max\{\frac{1}{2}, \frac{1}{2}\|Dg(\cdot, m)\|_{\infty, \R^{d}}\}$  and $g(x,\cdot)$ uniformly continuous w.r.t. space.
 \end{itemize}
 We also recall that $m_{0}:= \pi_{1} \sharp \mu_{0}$.  


Let $\Gamma$ be the set of $C^{1}$ curves $\gamma:[0,T] \to \R^{d}$, endowed with the local uniform convergence of the curve and its derivative, and given $(t,x,v) \in [0,T] \times \R^{2d}$ let $\Gamma_{t}(x,v)$ be the subset of $\Gamma$ such that $\gamma(t)=x$, $\dot\gamma(t)=v$. Similarly, let $\Gamma_{t}(x)$ be the subset of $\Gamma$ such that $\gamma(t)=x$. Define the functional $J^{\eps}_{t,T}: \Gamma \to \R$ 
 \begin{align*}
 J^{\eps}_{t,T}(\gamma)=& \int_{t}^{T}{\left(\frac{\eps}{2}|\ddot\gamma(s)|^{2} + L_{0}(\gamma(s), \dot\gamma(s), m^{\eps}_{s}) \right)\ ds} + g(\gamma(T),m^{\eps}_{T}), \quad \text{if}\,\, \gamma \in H^{2}(0,T; \R^{d})
 \end{align*}
 and set $J^{\eps}_{t,T}(\gamma)=+\infty$ if $\gamma \not\in H^{2}(0,T;\R^{d})$. 
 Then, from \cite{bib:CM, bib:YA} we know that there exist a solution $(u^{\eps}, \mu^{\eps}) \in W^{1,\infty}_{loc}([0,T] \times \R^{2d}) \times C([0,T]; \PP_{1}(\R^{2d}))$ to system \eqref{eq:epsMFG} such that
 \begin{equation}\label{eq:value}
 	u^{\eps}(t,x,v)= \inf_{\gamma \in \Gamma_{t}(x,v)} J^{\eps}_{t,T}(\gamma)
 \end{equation}
and for any $t \in [0,T]$ the probability measure $\mu^{\eps}_{t}$ is the image of $\mu_{0}$ under the flow
\begin{align}\label{eq:epsflow}
\begin{cases}
	\dot\gamma(t) =   v(t)
	\\
	\dot{v}(t) =  - \frac{1}{\eps}D_{v}u^{\eps}(t,\gamma(t), v(t)).
\end{cases}	
\end{align}
That is, $u^{\eps}$ solves the Hamilton-Jacobi equation in the viscosity sense and $\mu^{\eps}$ solves the continuity equation in the sense of distributions.  
\begin{remarks}\label{rem:mepsoptimal}\em
Note that for a.e. $(x,v) \in \R^{2d}$ there exists a unique solution to system \eqref{eq:epsflow}, which we will denote by $\gamma^{\eps}_{(x,v)}$,  such that $\gamma_{(x,v)}^{\eps}(0)=x$ and $\dot\gamma_{(x,v)}^{\eps}(0)=v$. Moreover, $\gamma^{\eps}_{(x,v)}(\cdot)$ is optimal for $u^{\eps}(t,x,v)$ and satisfies $\gamma_{(x,v)}^{\eps}(t)=x$, $\dot\gamma_{(x,v)}^{\eps}(t)=v$. 
\end{remarks}



\begin{theorem}[{\bf Main result 1}]\label{thm:main1}
	Assume {\bf (M1)} -- {\bf (M3)} and {\bf (BC1)}, {\bf (BC2)}. Let $(u^{\eps}, \mu^{\eps})$ be a solution to \eqref{eq:epsMFG} and let $m^{\eps}_{t}=\pi_{1} \sharp \mu^{\eps}_{t}$ for any $t \in [0,T]$. Then, there exists a sequence $\{\eps_{k}\}_{k \in \N}$ with $\eps_{k} \downarrow 0$, as $k \to \infty$, a function $u^{0} \in W^{1,\infty}_{loc}([0,T] \times \R^{d})$ and a flow of probability measures $\{m^{0}_{t}\}_{t \in [0,T]} \in C([0,T]; \PP_{1}(\R^{d}))$ such that for any $R \geq 0$
	\begin{equation*}
	\lim_{k \to \infty} u^{\eps_{k}}(t,x,v) = u^{0}(t,x), \quad \text{uniformly on}\,\, [0,T] \times \overline{B}_{R} \times \overline{B}_{R}
	\end{equation*}
	and 
	\[
	\lim_{k \to \infty} m^{\eps_{k}}=  m^{0}, \quad \text{in}\,\, C([0,T]; \PP_{1}(\R^{d})).
	\]
	Moreover, the following holds.
	\begin{itemize}
    \item[($i$)] $(u^{0},  m^{0}) \in W^{1,\infty}_{loc}([0,T] \times \R^{d}) \times C([0,T]; \PP_{1}(\R^{d}))$ is a solution of
	\begin{align}\label{eq:limitMFG}
	\begin{cases}
		-\partial_{t} u^{0}(t,x) + H_{0}(x, D_{x}u^{0}(t,x), m^{0}_{t})=0, & \quad (t,x) \in [0,T] \times \R^{d}
		\\
		\partial_{t}  m^{0}_{t} - \ddiv\Big( m^{0}_{t}D_{p}H_{0}(x, D_{x}u^{0}(t,x), m^{0}_{t}) \Big)=0, & \quad (t,x) \in [0,T] \times \R^{d}
		\\
		m^{0}_{0}= m_{0},\,\, u^{0}(T,x)=g(x,m^{0}_{T}), & \quad x \in \R^{d},
	\end{cases}
	\end{align}
that is, $u^{0}$ solves the Hamilton-Jacobi equation in the viscosity sense and $m^{0}$ is a solution of the continuity equation in the sense of distributions.
   \item[($ii$)]  For any $t \in [0,T]$ the probability measure $m^{0}_{t}$ is the image of $m_{0}$ under the Euler flow associated with $L_{0}$.  
	\end{itemize}
\end{theorem}

\begin{remarks}\em
Let $(u^{\eps}, \mu^{\eps})$ be a solution to \eqref{eq:epsMFG}. Assume that $H_{0}$ is of separated form, i.e., there exists a coupling function $F:\R^{d} \times \PP_{1}(\R^{d}) \to \R$ such that 
\begin{equation*}
	H_{0}(x,p,m)=H(x,p)-F(x,m), \quad \forall (x,p,m) \in \R^{2d} \times \PP_{1}(\R^{d}).
\end{equation*}
Moreover, assume that $F$ is continuous w.r.t. all variables, that the map $x \mapsto F(x,m)$ belongs to $C^{1}_{b}(\R^{d})$ and that the functions $F$, $g$ are monotone in the sense of Lasry-Lions, i.e. 
\begin{align*}
\int_{\R^{d}}{\big(F(x,m_{1})-F(x,m_{2}) \big)\ (m_{1}(dx)-m_{2}(dx))} \geq 0, & \quad  \forall\ m_{1}, m_{2} \in \PP_{1}(\R^{d})
\\
	\int_{\R^{d}}{\big(g(x,m_{1})-g(x,m_{2}) \big)\ (m_{1}(dx)-m_{2}(dx))} \geq 0, \quad &  \forall\ m_{1}, m_{2} \in \PP_{1}(\R^{d}).
\end{align*}
 Then, from \cite{bib:YA, bib:CM} we know that that there exists a unique solution $(u^{0},  m^{0}) \in W^{1,\infty}_{loc}([0,T] \times \R^{d}) \times C([0,T]; \PP_{1}(\R^{d}))$ of \eqref{eq:limitMFG} and thus as $(u^{\eps}, m^{\eps})$ is relatively compact  then convergence of $(u^{\eps}, m^{\eps})$ holds for the whole sequence. 
\end{remarks}

 \subsection{Convergence to MFG of control}

We now consider the function $L_{0}: \R^{d} \times \PP_{1}(\R^{2d}) \to \R$ and we assume the following. 
	 	\begin{itemize}
		\item[{\bf (C1)}] $L_{0}$ is continuous w.r.t. all variables and for any $\mu \in \PP_{1}(\R^{2d})$ the map $ x \mapsto L_{0}(x, \mu)$ belongs to $C^{1}(\R^{d})$.
       \item[{\bf (C2)}] There exist two moduli $\theta: \R_{+} \to \R_{+}$ and $\omega_{0}: \R_{+} \to \R_{+}$ such that
        \begin{equation*}
        |L_{0}(x,\mu_{1})-L_{0}(x,\mu_{2})| \leq \theta(|x|) \omega_{0}(d_{1}(\mu_{1},\mu_{2})),	
        \end{equation*}
for any $x \in \R^{d}$ and $\mu_{1}$, $\mu_{2} \in \PP_{1}(\R^{2d})$.
	\end{itemize}
	

We consider the MFG system  
	\begin{align}\label{eq:epsMFG}
 \begin{cases}
 	-\partial_{t} u^{\eps} +\frac{1}{2\eps}|D_{v}u^{\eps}|^{2} - \langle D_{x}u^{\eps}, v \rangle - \frac{1}{2}|v|^{2} - L_{0}(x,\mu^{\eps}_{t})= 0, & (t,x,v) \in [0,T] \times \R^{2d}
 	\\
 	\partial_{t}\mu^{\eps}_{t} - \langle D_{x}\mu^{\eps}_{t},v \rangle - \frac{1}{\eps}\ddiv_{v}\left(\mu^{\eps}_{t}D_{v}u^{\eps} \right)=0, &  (t,x,v) \in [0,T] \times \R^{2d}
\\
\mu^{\eps}_{0}=\mu_{0}, \quad u^{\eps}(T,x,v)=g(x,\mu^{\eps}_{T}), & (x,v) \in \R^{2d}
 \end{cases}	
 \end{align}
 
and we assume
 \begin{itemize}
 \item[{\bf (A1)}] the measure $\mu_{0} \in \PP(\R^{2d})$ is absolutely continuous w.r.t. Lebesgue measure, we still denote by $\mu_{0}$ its density, and it has compact support. 
 \item[{\bf (A2)}] The terminal costs $g(\cdot, \mu)$ belongs to $C^{1}_{b}(\R^{d})$, $g(x,\cdot)$ uniformly continuous w.r.t. space and we have that $M_{0} \geq \max\{\frac{1}{2}, \frac{1}{2}\|Dg(\cdot, \mu)\|_{\infty, \R^{d}}\}$ .
 \end{itemize}



Similarly to the previous part, we define the functional $J^{\eps}_{t,T}: \Gamma \to \R$ 
 \begin{align*}
 J^{\eps}_{t,T}(\gamma)=& \int_{t}^{T}{\left(\frac{\eps}{2}|\ddot\gamma(s)|^{2} + \frac{1}{2}|\dot\gamma(s)|^{2}+L_{0}(\gamma(s), \mu^{\eps}_{s}) \right)\ ds} + g(\gamma(T),\mu^{\eps}_{T}), \quad \text{if}\,\, \gamma \in H^{2}(0,T; \R^{d})
 \end{align*}
 and set $J^{\eps}_{t,T}(\gamma)=+\infty$ if $\gamma \not\in H^{2}(0,T;\R^{d})$. 
 Then, from \cite{bib:CM, bib:YA} we know that there exist a solution $(u^{\eps}, \mu^{\eps}) \in W^{1,\infty}_{loc}([0,T] \times \R^{2d}) \times C([0,T]; \PP_{1}(\R^{2d}))$ to system \eqref{eq:epsMFG} such that
 \begin{equation}\label{eq:value}
 	u^{\eps}(t,x,v)= \inf_{\gamma \in \Gamma_{t}(x,v)} J^{\eps}_{t,T}(\gamma)
 \end{equation}
and for any $t \in [0,T]$ the probability measure $\mu^{\eps}_{t}$ is the image of $\mu_{0}$ under the flow
\begin{align}\label{eq:epsflow2}
\begin{cases}
	\dot\gamma(t) =   v(t)
	\\
	\dot{v}(t) =  - \frac{1}{\eps}D_{v}u^{\eps}(t,\gamma(t), v(t)).
\end{cases}	
\end{align}
That is, $u^{\eps}$ solves the Hamilton-Jacobi equation in the viscosity sense and $\mu^{\eps}$ solves the continuity equation in the sense of distributions.  



\begin{theorem}[{\bf Main result 2}]\label{thm:main2}
	Assume {\bf (C1)} -- {\bf (C3)} and {\bf (A1)}, {\bf (A2)}. Let $(u^{\eps}, \mu^{\eps})$ be a solution to \eqref{eq:epsMFG}. Then, there exists a sequence $\{\eps_{k}\}_{k \in \N}$ with $\eps_{k} \downarrow 0$, as $k \to \infty$, a function $u^{0} \in W^{1,\infty}_{loc}([0,T] \times \R^{d})$ and a flow of probability measures $\{\mu^{0}_{t}\}_{t \in [0,T]} \in C([0,T]; \PP_{1}(\R^{2d}))$ such that for any $R \geq 0$
	\begin{equation*}
	\lim_{k \to \infty} u^{\eps_{k}}(t,x,v) = u^{0}(t,x), \quad \text{uniformly on}\,\, [0,T] \times \overline{B}_{R} \times \overline{B}_{R}
	\end{equation*}
	and 
	\[
	\lim_{k \to \infty} \mu^{\eps_{k}}=  \mu^{0}, \quad \text{in}\,\, C([0,T]; \PP_{1}(\R^{2d})).
	\]
	Moreover, we have that the pair $(u^{0},  \mu^{0}) \in W^{1,\infty}_{loc}([0,T] \times \R^{d}) \times C([0,T]; \PP_{1}(\R^{d}))$ is a solution of the MFG of control system
	\begin{align}\label{eq:limitMFG}
	\begin{cases}
		(i)\,\, -\partial_{t} u^{0}(t,x) + \frac{1}{2}|D_{x}u^{0}(t,x)|^{2}-  L_{0}(x, \mu^{0}_{t})=0, & \quad (t,x) \in [0,T] \times \R^{d}
		\\
		(ii)\,\, \partial_{t}  m^{0}_{t} - \ddiv\big( m^{0}_{t}D_{x}u^{0}(t,x) \big)=0, & \quad (t,x) \in [0,T] \times \R^{d}
		\\
		(iii)\,\, \mu^{0}_{t} = (\text{Id}(\cdot), D_{x}u^{0}(t, \cdot)) \sharp m^{0}_{t}, & \quad t \in [0,T]
		\\
		m^{0}_{0}= m_{0},\,\, u^{0}(T,x)=g(x,\mu^{0}_{T}), & \quad x \in \R^{d},
	\end{cases}
	\end{align}
that is, $u^{0}$ solves the Hamilton-Jacobi equation in the viscosity sense and $m^{0}_{t}=\pi_{1} \sharp \mu^{0}_{t}$, for all $t \in [0,T]$, is a solution of the continuity equation in the sense of distributions. Furthermore, the whole measure $\mu^{0}$ is given by ($iii$) in \eqref{eq:limitMFG}.
\end{theorem}

\section{Proof of \Cref{thm:main1}} 
\label{sec:mainresult}

	In order to prove \Cref{thm:main1} we proceed by steps analyzing the behavior of the value function $u^{\eps}$ and that of the flow of probability measures $\{m^{\eps}_{t}\}_{t \in [0,T]}$ separately. First, we show that $u^{\eps}$ is equibounded and we prove that, up to a subsequence, $m^{\eps}$ converges to a flow of probability measure in $C([0,T]; \PP_{1}(\R^{d}))$. Then, we address the convergence of the value function, up to a subsequence, to a solution of a suitable Hamilton-Jacobi equation and we study the limit of its minimizing trajectories. Finally, we are able to characterize the limit flow of measures as solution of a continuity equation which coupled with the Hamilton-Jacobi equation, previously constructed, define the limit MFG system \eqref{eq:limitMFG}.

		\begin{lemma}\label{lem:Mequibound}
Assume {\bf (M1)} -- {\bf (M3)} and {\bf (BC1)}, {\bf (BC2)}. Then we have that 
\begin{equation*}
	-TM_{0}-\| g(\cdot, m^{\eps}_{T})\|_{\infty, \R^{d}} \leq u^{\eps}(t,x,v) \leq M_{0}T(1+|v|^{2}) + \| g(\cdot, m^{\eps}_{T})\|_{\infty, \R^{d}}, 
\end{equation*}
for any $(t,x,v) \in [0,T] \times \R^{2d}$ and for any $\eps > 0$.
\end{lemma}
\proof
First, since $u^{\eps}$ satisfy \eqref{eq:value}, from \eqref{eq:Mcontrolgrowth} and {\bf (BC)} follows that for any $(t,x,v) \in [0,T] \times \R^{2d}$ there holds
\begin{equation*}
u^{\eps}(t,x,v) \geq - C_{0}T-\|g(\cdot, m^{\eps}_{T})\|_{\infty, \R^{d}}. 	
\end{equation*}
On the other hand, let us recall that $u^{\eps}$ solves the Hamilton-Jacobi equation 
\begin{equation}\label{eq:uhj}
-\partial_{t} u^{\eps} +\frac{1}{2\eps}|D_{v}u^{\eps}|^{2} - \langle D_{x}u^{\eps}, v \rangle -L_{0}(x, v, m^{\eps}_{t})= 0. 
\end{equation}
Then, the function
\begin{equation*}
\zeta(t,x,v)= g(x,m^{\eps}_{T}) +C(1+|v|^{2})(T-t), \quad (t,x,v) \in [0,T] \times \R^{2d}	
\end{equation*}
is a supersolution to \eqref{eq:uhj} for a suitable choice of the real constant $C \geq 0$. Indeed, we have that \begin{align*}
- \ &\partial_{t}\zeta(t,x,v) + \frac{1}{2\eps}|D_{v}\zeta(t,x,v)|^{2} - \langle D_{x}\zeta(t,x,v), v \rangle -L_{0}(x,v)
\\
\geq\ &    C(1+|v|^{2}) + 2\frac{(T-t)^{2}C^{2}}{\eps}|v|^{2} -\langle D_{x}g(x,m^{\eps}_{T}), v\rangle - M_{0}(1+|v|^{2})
\\
\geq\ &  C(1+|v|^{2}) -\frac{1}{2}\|Dg(\cdot,m^{\eps}_{T})\|_{\infty,\R^{d}} - \frac{1}{2}|v|^{2} - M_{0}(1+|v|^{2})
\end{align*}
where the last inequality holds by Young's inequality. 
Thus, taking $C=2M_{0}$ by {\bf (BC)}  we obtain  
\begin{equation*}
	M_{0}(1+|v|^{2}) -\frac{1}{2}\|Dg(\cdot,m^{\eps}_{T})\|_{\infty,\R^{d}} - \frac{1}{2}|v|^{2}  \geq 0.
\end{equation*}
So, we get the result by Comparison Theorem \cite[Theorem 2.12]{bib:BD}. \qed

		
		\begin{corollary}\label{cor:Mapproxbuondedvel}
	Assume {\bf (M1)} -- {\bf (M3)} and {\bf (BC1)}, {\bf (BC2)}. Let $(t,x,v) \in [0,T] \times \R^{2d}$ and let $\gamma^{\eps}$ be a minimizer for $u^{\eps}(t,x,v)$. Then, there exists a constant $Q_{1} \geq 0$ such that 
	\begin{equation*}
	\int_{t}^{T}{|\dot\gamma^{\eps}(s)|^{2}\ ds} \leq Q_{1}(1+|v|^{2}), \quad \forall\ \eps > 0.
	\end{equation*}
where $Q_{1}$ is independent of $\eps$, $t$, $x$ and $v$. 
\end{corollary}
\proof
On the one hand, from \Cref{lem:Mequibound} we know that 
\begin{equation*}
u^{\eps}(t,x,v) \leq {M}_{0}T\big(1+|v|^{2}\big) + \|g(\cdot, m^{\eps}_{T})\|_{\infty, \R^{d}}, \quad \forall\ (t,x,v) \in [0,T] \times \R^{2d}.	
\end{equation*}
On the other hand, let $(t,x,v) \in [0,T] \times \R^{2d}$ and let $\gamma^{\eps}$ be a minimizer for $u^{\eps}(t,x,v)$. Then, by \eqref{eq:Mcontrolgrowth} we have that 
\begin{align*}
u^{\eps}(t,x,v)=\ &\int_{t}^{T}{\left(\frac{\eps}{2}|\ddot\gamma^{\eps}(s)|^{2} + L_{0}(\gamma^{\eps}(s), \dot\gamma^{\eps}(s), m^{\eps}_{s}) \right)\ ds} + g(\gamma^{\eps}(T),m^{\eps}_{T})
\\
\geq & \int_{t}^{T}{L_{0}(\gamma^{\eps}(s), \dot\gamma^{\eps}(s), m^{\eps}_{s})\ ds} -\|g(\cdot, m^{\eps}_{T})\|_{\infty, \R^{d}} 
\\
\geq\ & \int_{t}^{T}{\left(\frac{1}{M_{0}}|\dot\gamma^{\eps}(s)|^{2} - M_{0} \right)\ ds} -\|g(\cdot, m^{\eps}_{T})\|_{\infty, \R^{d}}.	
\end{align*}
Therefore, combining the above inequalities we get
\begin{equation*}
\int_{t}^{T}{|\dot\gamma^{\eps}(s)|^{2}\ ds}	 \leq  2M_{0}\big(\|g(\cdot, m^{\eps}_{T})\|_{\infty, \R^{d}} + M_{0}T(1+|v|^{2})\big)=:Q_{1}(1+|v|^{2})
\end{equation*}
where $Q_{1}$ depends only on $M_{0}$, $T$ and  $\|g(\cdot, m^{\eps}_{T})\|_{\infty, \R^{d}}$ which is bounded uniformly in $m^{\eps}_{T}$. \qed

\begin{corollary}\label{cor:Holder}
	Assume {\bf (M1)} -- {\bf (M3)} and {\bf (BC1)}, {\bf (BC2)}. Then, there exists a constant $Q_{2} \geq 0$ such that for any $s_{1}$, $s_{2} \in [0,T]$ with $s_{1} \leq s_{2}$ there holds
	\begin{equation*}
	d_{1}(m^{\eps}_{s_{2}}, m^{\eps}_{s_{1}}) \leq Q_{2}|s_{1}-s_{2}|^{\frac{1}{2}}, \quad \forall\ \eps > 0
	\end{equation*}
	where $Q_{2}$ is independent of $\eps$. 
\end{corollary}
\proof
We first recall that for any $t \in [0,T]$ we know that $m^{\eps}_{t}=\pi_{1} \sharp \mu^{\eps}_{t}$ where $\mu^{\eps}_{t}$ is the image of $\mu_{0}$ under the flow \eqref{eq:epsflow} whose space marginal we denote by $\gamma^{\eps}_{(x,v)}$ for $(x,v) \in \R^{2d}$. 

Let $s_{1}$, $s_{2} \in [0,T]$ be such that $s_{1} \leq s_{2}$. Then, by \eqref{eq:100} we have that 
\begin{align*}
d_{1}(m^{\eps}_{s_{1}}, m^{\eps}_{s_{2}}) \leq \int_{\R^{d}}{|\gamma^{\eps}_{(x,v)}(s_{1})-\gamma^{\eps}_{(x,v)}(s_{2})|\ \mu_{0}(dx,dv)}
\end{align*}
and thus, appealing to \Cref{cor:Mapproxbuondedvel} and the H\"older inequality we obtain 
\begin{align*}
	d_{1}(m^{\eps}_{s_{1}}, m^{\eps}_{s_{2}}) \leq |s_{1}-s_{2}|^{\frac{1}{2}}\left(\int_{\R^{2d}}{Q_{1}(1+|v|^{2})\ \mu_{0}(dx,dv)}\right)^{\frac{1}{2}}.
\end{align*}
So, since  $\mu_{0}$ has compact support we get the result setting 
\begin{equation*}
	Q_{2}=\left(\int_{\R^{2d}}{Q_{1}(1+|v|^{2})\ \mu_{0}(dx,dv)}\right)^{\frac{1}{2}}. \eqno\square
\end{equation*}

We are now ready to prove that the flow of probability measures $m^{\eps}$ converges, up to a subsequence. First, we recall that for any $t \in [0,T]$ the measure $m_{t}^{\eps}$ is the space marginal of $\mu^{\eps}_{t}$ which is given by the push-forward of the initial distribution $\mu_{0}$ by the optimal flow \eqref{eq:epsflow}, that is 
\begin{align*}
\begin{cases}
	\dot\gamma_{(x,v)}(t) =   v(t), &  \gamma_{(x,v)}(0)=x
	\\
	\dot{v}(t) =  - \frac{1}{\eps}D_{v}u^{\eps}(t,\gamma_{(x,v)}(t), v(t)), &  v(0)=v.
\end{cases}	
\end{align*}

\begin{theorem}\label{thm:spacetightness}
			Assume {\bf (M1)} -- {\bf (M3)} and {\bf (BC1)}, {\bf (BC2)}. Then, the flow of measures $\{m^{\eps}_{t}\}_{t \in [0,T]}$ is tight and there exists a sequence $\{\eps_{k}\}_{k \in \N}$ such that $m^{\eps_{k}}$ converges to some probability measure $m^{0}$ in $C([0,T]; \PP_{1}(\R^{d}))$. 
		\end{theorem}
		\proof
		Since $m_{t}^{\eps}=\pi_{1} \sharp \mu^{\eps}_{t}$, for any $t \in [0,T]$,  where $\mu^{\eps}_{t} $ is given by push-forward of $\mu_{0}$ under the flow \eqref{eq:epsflow}, we know that 
		\begin{equation*}
		\int_{\R^{d}}{|x|^{2}\ m_{t}^{\eps}(dx)}=\int_{\R^{2d}}{|\gamma_{(x,v)}^{\eps}(t)|^{2}\ \mu_{0}(dx,dv)}.
		\end{equation*}
		So, we are interested in estimating the curve $\gamma_{(x,v)}^{\eps}$ for any $(x,v)$, uniformly in $\eps >0$. In order to get it, from \Cref{cor:Mapproxbuondedvel} we immediately deduce that 
		\begin{equation*}
		|\gamma^{\eps}_{(x,v)}(s)| \leq |x|	+ \sqrt{T}\sqrt{Q_{1}}(1+|v|^{2})^{\frac{1}{2}}, \quad \forall\ s \in [0,T].
		\end{equation*}
	 Hence, for any $t \geq 0$ we have that
		\begin{align*}
		& \int_{\R^{d}}{|x|^{2}\ m^{\eps}_{t}(dx)} = \int_{\R^{d}}{|\gamma^{\eps}_{(x,v)}(t)|^{2}\  \mu_{0}(dx,dv)}
		\\
		\leq & \int_{\R^{2d}}{C_{0}\big(|x|^{2} + TQ_{1}(1+|v|^{2})\big)\ \mu_{0}(dx,dv)}	
		\end{align*}
		for some constant $C_{0} \geq 0$. Thus, since $\mu_{0}$ has compact support we deduce that $\{m^{\eps}_{t}\}_{t \in [0,T]}$ has bounded second-order momentum, uniformly in $\eps > 0$ and, consequently, $\{m^{\eps}_{t}\}_{t \in [0,T]}$ is tight. Therefore, by Prokhorov Theorem and Ascoli-Arzela Theorem, $\mu_{0}$ has uniformly bounded support and by \Cref{cor:Holder} $m^{\eps}_{t}$ is equicontinuous in time, there exists a sequence $\{\eps_{k}\}_{k \in \N}$ and measure $m^{0} \in C([0,T]; \PP_{1}(\R^{d}))$ such that $m^{\eps_{k}} \to m^{0}$ in $C([0,T]; \PP_{1}(\R^{d}))$.\qed

\smallskip
Next, we turn to the convergence of the value function $u^{\eps}$. Before proving it, we need preliminary estimates on the oscillation of the value function w.r.t. velocity variable and then w.r.t. time and space variable. In particular, we will show that the function $u^{\eps}(t,x,\cdot)$ has decreasing oscillation w.r.t. $\eps$, which will allowed us to conclude that the limit function does not depend on $v$.

\begin{lemma}\label{lem:controllability}
	Assume {\bf (M1)} -- {\bf (M3)} and {\bf (BC1)}, {\bf (BC2)}. Let $R \geq 0$ and let $(x,v_{0})$, $(x,v) \in \R^{d} \times \overline{B}_{R}$. Then, there exists $C_{R} \geq 0$ and a parametric curve $\sigma: [0,\sqrt{\eps}] \to \R^{d}$ such that 
	\begin{equation*}
		\sigma(0)=\sigma(\sqrt{\eps}) = x, \quad \dot\sigma(0)=v_{0}, \quad \dot\sigma(\sqrt{\eps})=v
	\end{equation*}
and
	\begin{equation*}
	\frac{1}{\sqrt{\eps}}\int_{0}^{\sqrt{\eps}}{\left(\frac{\eps}{2}|\ddot\sigma(s)|^{2} + L_{0}(\sigma(s), \dot\sigma(s), m^{\eps}_{s}) \right)\ ds} \leq C_{R}
	\end{equation*}
where $C_{R}$ is independent of $\eps$, $x$, $v$ and $v_{0}$. 
\end{lemma}
\proof
Let $R \geq 0$ and let $(x,v_{0})$, $(x,v) \in \R^{d} \times \overline{B}_{R}$. Define the curve $\sigma:[0,\sqrt{\eps}] \to \R^{d}$ by 
\begin{equation*}
\sigma(t)=x+v_{0}t+Bt^{2} + At^{3}	
\end{equation*}
with $A$, $B \in \R$ satisfying the following conditions
\begin{equation*}
	\sigma(0)=\sigma(\sqrt{\eps}) = x, \quad \dot\sigma(0)=v_{0}, \quad \dot\sigma(\sqrt{\eps})=v.
\end{equation*}
Thus, we obtain 
\begin{equation*}
\begin{cases}
	B= -(2v_{0}+v)\eps^{-\frac{1}{2}},
	\\
	A= (v+v_{0})\eps^{-1}.
\end{cases}	
\end{equation*}
Hence, we get
\begin{multline*}
 \int_{0}^{\sqrt{\eps}}{\left(\frac{\eps}{2}|\ddot\sigma(s)|^{2} + L_{0}(\sigma(s), \dot\sigma(s), m^{\eps}_{s})\right)\ ds} \\ \leq \int_{0}^{\sqrt{\eps}}{\left(\frac{\eps}{2}|2B+6At|^{2} + M_{0}(1+|v+2tB+3t^{2}A|^{2})\right)\ ds} \leq \widehat{C}\sqrt{\eps}R^{2}
\end{multline*}
for some positive constant $\widehat{C}$ and the proof is thus complete. \qed

		\begin{lemma}\label{lem:oscillationvel}
Assume {\bf (M1)}, {\bf (M2)} and {\bf (BC)}. Let $R \geq 0$, let $T > 1$ and $\eps > 0$. Then, there exists $\widehat{C}_{R}(\eps) \geq 0$ such that for any $t \in [0,T]$, any $x \in \R^{d}$, and any $v$, $w$ in $\overline{B}_{R}$  there holds
\begin{equation*}
|u^{\eps}(t,x,v)-u^{\eps}(t,x,w)| \leq \widehat{C}_{R}(\eps)	
\end{equation*}
and $\widehat{C}_{R}(\eps) \to 0$ as $\eps \downarrow 0$. 
\end{lemma}
\proof
Fix $R \geq 0$ and take $(x,v)$, $(x,w) \in \R^{d} \times \overline{B}_{R}$. Let $\gamma^{\eps}$ be a minimizer for $u^{\eps}(t,x,v)$ and define the curve
\begin{align*}
\widehat\gamma(s)=
\begin{cases}
	\sigma(s-t), & \quad s \in [t,t+\sqrt{\eps}]
	\\
	\gamma^{\eps}(s-\sqrt{\eps}), & \quad s \in [t+\sqrt{\eps},T]
\end{cases}	
\end{align*}
where $\sigma: [0,\sqrt{\eps}] \to \R^{2d}$ connects, in the sense of  \Cref{lem:controllability}, $(x,w)$ with $(x,v)$. Then, we obtain
\begin{multline*}
u^{\eps}(t,x,w) - u^{\eps}(t,x,v) \leq\  \int_{t}^{t+\sqrt{\eps}}{\left(\frac{\eps}{2}|\ddot\sigma(s-t)|^{2} + L_{0}(\sigma(s-t), \dot\sigma(s-t), m^{\eps}_{s}) \right)\ ds}  
\\
+\  \int_{t+\sqrt{\eps}}^{T}{\left(\frac{\eps}{2}|\ddot\gamma^{\eps}(s-\sqrt{\eps})|^{2} + L_{0}(\gamma^{\eps}(s-\sqrt{\eps}), \dot\gamma^{\eps}(s-\sqrt{\eps}), m^{\eps}_{s}) \right)\ ds} 
\\
+\  g(\gamma^{\eps}(T-\sqrt{\eps}), m^{\eps}_{T}) - u^{\eps}(t,x,v)
\\
=\  \int_{t}^{t+\sqrt{\eps}}{\left(\frac{\eps}{2}|\ddot\sigma(s-t)|^{2} + L_{0}(\sigma(s-t), \dot\sigma(s-t), m^{\eps}_{s-t}) \right)\ ds} 
\\
+\  \int_{t+\sqrt{\eps}}^{T}{\left(\frac{\eps}{2}|\ddot\gamma^{\eps}(s-\sqrt{\eps})|^{2} + L_{0}(\gamma^{\eps}(s-\sqrt{\eps}), \dot\gamma^{\eps}(s-\sqrt{\eps}), m^{\eps}_{s-\sqrt{\eps}}) \right)\ ds} 
\\
+\  g(\gamma^{\eps}(T), m^{\eps}_{T}) + g(\gamma^{\eps}(T-\sqrt{\eps}), m^{\eps}_{T})-g(\gamma^{\eps}(T), m^{\eps}_{T}) - u^{\eps}(t,x,v)
\\
+\   \int_{t}^{t+\sqrt{\eps}}{\big(L_{0}(\sigma(s-t), \dot\sigma(s-t), m^{\eps}_{s}) -L_{0}(\sigma(s-t), \dot\sigma(s-t), m^{\eps}_{s-t})\big)\ ds}
\\
+\  \int_{t+\sqrt{\eps}}^{T}{\big(L_{0}(\gamma^{\eps}(s-\sqrt{\eps}), \dot\gamma^{\eps}(s-\sqrt{\eps}), m^{\eps}_{s}) - L_{0}(\gamma^{\eps}(s-\sqrt{\eps}), \dot\gamma^{\eps}(s-\sqrt{\eps}), m^{\eps}_{s-\sqrt{\eps}} \big)\ ds}.
\end{multline*}
Now, from \Cref{lem:controllability} we know that 
\begin{equation}\label{eq:111}
	\int_{t}^{t+\sqrt{\eps}}{\left(\frac{\eps}{2}|\ddot\sigma(s-t)|^{2} + L_{0}(\sigma(s-t), \dot\sigma(s-t), m^{\eps}_{s-t}) \right)\ ds} \leq C_{R}\sqrt{\eps},
\end{equation}
and, moreover, from the optimality of $\gamma^{\eps}$ we get 
\begin{align}\label{eq:112}
\begin{split}
	& \int_{t+\sqrt{\eps}}^{T}{\left(\frac{\eps}{2}|\ddot\gamma^{\eps}(s-\sqrt{\eps})|^{2} + L_{0}(\gamma^{\eps}(s-\sqrt{\eps}), \dot\gamma^{\eps}(s-\sqrt{\eps}), m^{\eps}_{s-\sqrt{\eps}}) \right)\ ds} -u^{\eps}(t,x,v)
	\\
	\leq &  - \int_{T-\sqrt{\eps}}^{T}{\left(\frac{\eps}{2}|\ddot\gamma^{\eps}(s)|^{2}+L_{0}(\gamma^{\eps}(s), \dot\gamma^{\eps}(s), m^{\eps}_{s}) \right)\ ds} \leq 0.
	\end{split}
\end{align}
Then, as observed before from \Cref{cor:Mapproxbuondedvel} we have that 
\begin{equation*}
	|\gamma^{\eps}(s)| \leq |x| + \sqrt{T}\sqrt{Q_{1}}(1+|v|^{2})^{\frac{1}{2}}, \quad \forall\ s \in [0,T]
\end{equation*}
and also that the curve $\sigma$ is bounded. 
 Hence, by {\bf (M3)} and \Cref{cor:Holder} we deduce that there exists $P(\eps) \geq 0$, with $P(\eps) \to 0$ as $\eps \downarrow 0$, such that 
\begin{align}\label{eq:113}
\begin{split}
	&  \int_{t}^{t+\sqrt{\eps}}{\big(L_{0}(\sigma(s-t), \dot\sigma(s-t), m^{\eps}_{s}) -L_{0}(\sigma(s-t), \dot\sigma(s-t), m^{\eps}_{s-t})\big)\ ds}
\\
+\ & \int_{t+\sqrt{\eps}}^{T}{\big(L_{0}(\gamma^{\eps}(s-\sqrt{\eps}), \dot\gamma^{\eps}(s-\sqrt{\eps}), m^{\eps}_{s}) - L_{0}(\gamma^{\eps}(s-\sqrt{\eps}), \dot\gamma^{\eps}(s-\sqrt{\eps}), m^{\eps}_{s-\sqrt{\eps}} \big)\ ds}
\\
+\ & g(\gamma^{\eps}(T-\sqrt{\eps}), m^{\eps}_{T})-g(\gamma^{\eps}(T), m^{\eps}_{T}) \leq P(\eps)
\end{split}
\end{align}
where we have used that the modulus $\theta$ in {\bf (M3)} is bounded from the boundedness of $\gamma^{\eps}$ and $\sigma$. 
Therefore, combining \eqref{eq:111}, \eqref{eq:112} and \eqref{eq:113} we get the result. \qed

\begin{proposition}\label{prop:Mequicontinuity}
Assume {\bf (M1)} -- {\bf (M3)} and {\bf (BC1)}, {\bf (BC2)}. Then, for any $R \geq 0$ there exists a modulus $\omega_{R}:\R_{+} \to \R_{+}$ and a constant $C_{1} \geq 0$, independent of $R$, such that for any $\eps > 0$ the following holds:
\begin{align}
|u^{\eps}(t,x,v) - u^{\eps}(s,x,v)| \leq\ & \omega_{R}(|t-s|), \quad \forall\ (t,s,x,v) \in [0,T] \times [0,T] \times \overline{B}_{R} \times \overline{B}_{R}\label{eq:Mtimegradient}
	\\
|D_{x}u^{\eps}(t,x,v)| \leq\ & C_{1} T(1+|v|^{2}), \quad \text{a.e.}\,\, (t,x,v) \in [0,T] \times \R^{d} \times \R^{d}.\label{eq:Mspacegradient}	
\end{align} 
\end{proposition}
\proof
We begin by proving \eqref{eq:Mspacegradient}. Let $(t,x,v) \in [0,T] \times \R^{d} \times \R^{d}$ and let $\gamma^{\eps}$ be a minimizer for $u^{\eps}(t,x,v)$. Then, from \eqref{eq:MLspacegradient} we get
\begin{align*}
u^{\eps}(t,x+h,v) \leq & \int_{t}^{T}{\left(\frac{\eps}{2}|\ddot\gamma^{\eps}(s)|^{2} + L_{0}(\gamma^{\eps}(s)+h, \dot\gamma^{\eps}(s), m^{\eps}_{s}) \right)\ ds} + g(\gamma^{\eps}(T)+h, m^{\eps}_{T})
\\
=\ & u^{\eps}(t,x,v) + \int_{t}^{T}{\big(L_{0}(\gamma^{\eps}(s)+h, \dot\gamma^{\eps}(s), m^{\eps}_{s})-L_{0}(\gamma^{\eps}(s), \dot\gamma^{\eps}(s), m^{\eps}_{s})\big)\ ds}
\\
+\ & g(\gamma^{\eps}(T)+h, m^{\eps}_{T})- g(\gamma^{\eps}(T), m^{\eps}_{T})
\\
\leq\ & u^{\eps}(t,x,v) + \int_{t}^{T}{M_{0}|h|(1+|\dot\gamma^{\eps}(s)|^{2})\ ds} + \|Dg(\cdot, m^{\eps}_{T})\|_{\infty, \R^{d}}|h|.
\end{align*}
Hence, \Cref{cor:Mapproxbuondedvel} yields to the conclusion. 

 Next, we proceed to show \eqref{eq:Mtimegradient}. Let $R \geq 0$ and take $(t,x,v) \in [0,T] \times \overline{B}_{R} \times \overline{B}_{R}$. Let $\gamma^{\eps}$ be a minimizer for $u^{\eps}(t,x,v)$ and let $h \in [0,T-t]$. Then, we have that 
\begin{align*}
& u^{\eps}(t+h,x,v) 
\\
\leq\ & \int_{t+h}^{T}{\left(\frac{\eps}{2}|\ddot\gamma^{\eps}(s-h)|^{2} + L_{0}(\gamma^{\eps}(s-h), \dot\gamma^{\eps}(s-h), m^{\eps}_{s}) \right)\ ds} + g(\gamma(T-h), m^{\eps}_{T})
\\
=\ & \int_{t+h}^{T}{\left(\frac{\eps}{2}|\ddot\gamma^{\eps}(s-h)|^{2} + L_{0}(\gamma^{\eps}(s-h), \dot\gamma^{\eps}(s-h), m^{\eps}_{s-h}) \right)\ ds} + g(\gamma^{\eps}(T), m^{\eps}_{T})
\\
+\ & \int_{t+h}^{T}{\big(L_{0}(\gamma^{\eps}(s-h), \dot\gamma^{\eps}(s-h), m^{\eps}_{s}) -  L_{0}(\gamma^{\eps}(s-h), \dot\gamma^{\eps}(s-h), m^{\eps}_{s-h}) \big)\ ds} 
\\
+\ & g(\gamma^{\eps}(T-h), m^{\eps}_{T})-g(\gamma^{\eps}(T), m^{\eps}_{T})
\\
\leq\ & u^{\eps}(t,x,v) + \int_{t+h}^{T}{\theta(|\gamma^{\eps}(s-h)|)\omega_{0}(d_{1}(m^{\eps}_{s}, m^{\eps}_{s-h}))\ ds} + \|Dg(\cdot, m^{\eps}_{T})\|_{\infty, \R^{d}}|h|
\end{align*}
where the last inequality holds by {\bf (M3)}. Hence, from \Cref{cor:Mapproxbuondedvel} we know that 
\begin{equation*}
	|\gamma^{\eps}(s)| \leq |x| + \sqrt{T}\sqrt{Q_{1}}(1+|v|^{2})^{\frac{1}{2}}, \quad \forall\ s \in [0,T]
\end{equation*}
and thus $\theta(\cdot)$ turns out to be bounded. Therefore, appealing to \Cref{cor:Holder} we obtain
\begin{equation}\label{eq:Mone}
u^{\eps}(t+h,x,v) - u^{\eps}(t,x,v) \leq T\theta(R)\omega_{0}(|h|^{\frac{1}{2}}) +	\|Dg(\cdot, m^{\eps}_{T})\|_{\infty, \R^{d}}|h|. 
\end{equation}
On the other hand, let $R \geq 0$ and let $(t,x,v) \in [0,T] \times \overline{B}_{R}\times \overline{B}_{R}$. For $h \in [0,T-t]$, define the curve $\gamma:[t,t+h] \to \R^{d}$ by $\gamma(s)=x+(s-t)v$. Then, by Dynamic Programming Principle we deduce that 
\begin{align}\label{eq:Mtwo}
\begin{split}
 u^{\eps}(t,x,v) \leq\ & \int_{t}^{t+h}{\left(\frac{\eps}{2}|\ddot\gamma(s)|^{2} + L_{0}(\gamma(s), \dot\gamma(s), m^{\eps}_{s}) \right)\ ds}	 
\\
+\ & u^{\eps}(t+h, \gamma(t+h), \dot\gamma(t+h))
\\
=\ & \int_{t}^{t+h}{L_{0}(x+(s-t)v, v, m^{\eps}_{s})\ ds}	 + u^{\eps}(t+h, x+hv, v)
\\
\leq\ & M_{0}(1+R^{2})|h| + u^{\eps}(t+h,x,v) + C_{1}T(1+R^{2})|h|
\end{split}
\end{align}
 where we applied \eqref{eq:Mcontrolgrowth} and   \Cref{eq:Mspacegradient} to get the last inequality. Therefore, combining \eqref{eq:Mone} and \eqref{eq:Mtwo} the proof is complete. \qed

	\begin{remarks}\label{rem:modulustheta}\em
		Next, we study the behavior of the value function $u^{\eps}$ as $\eps \to 0$ and before doing that we recall the following argument needed to get uniform convergence to a function which does not depend on $v$. Assume that there exists a nonnegative function $\Theta(\delta_{0}, \eps_{0}, R_{0})$ such that 
		\begin{equation*}
			\Theta(\delta_{0}, \varepsilon_{0}, R_{0}) \to 0, \quad \text{as}\quad \varepsilon_{0}, \delta_{0} \downarrow 0,
		\end{equation*}
and assume that for any $|t_{1}-t_{2}| + |x_{1} - x_{2}| \leq \delta_{0}$, any $\eps \leq \eps_{0}$ and any  $|x_{i}|$, $|v_{i}| \leq R_{0}$  ($i=1,2$) there holds
		\begin{equation*}
		|u^{\eps}(t_{1}, x_{1}, v_{1}) - u^{\eps}(t_{2}, x_{2}, v_{2})| \leq \Theta(\delta_{0}, \varepsilon_{0}, R_{0}).
		\end{equation*}
	Then: if $u^{\eps}$ converge point-wise then $u^{\eps}$ converges locally uniformly and the limit function does not depend on $v$. \qed
		\end{remarks}

		Let $m^{0} \in C([0,T]; \PP_{1}(\R^{d}))$ be the flow of measures obtained in \Cref{thm:spacetightness} as limit of the flow $m^{\eps_{k}}$ in $C([0,T]; \PP_{1}(\R^{d}))$ for some subsequence $\eps_{k} \downarrow 0$. Define the function $u^{0}: [0,T] \times \R^{2d} \to \R$ by 
\begin{equation}\label{eq:Mlimitvaluefunction}
	u^{0}(t,x)= \inf_{\gamma \in \Gamma_{t}(x)} \left\{\int_{t}^{T}{L_{0}(\gamma(s), \dot\gamma(s), m^{0}_{s})\ ds} + g(\gamma(T), m^{0}_{T})\right\}.
	\end{equation}
		We will prove now that for the subsequence $\eps_{k}$ the sequence of value functions $u^{\eps_{k}}$ locally uniformly converge to $u^{0}$.

\begin{theorem}\label{thm:Meps}
Assume {\bf (M1)} -- {\bf (M3)} and {\bf (BC1)}, {\bf (BC2)}. Then, there exists a subsequence $\eps_{k} \downarrow 0$ such that  
	$u^{\eps_{k}}$ locally uniformly converges to $u^{0}$.
\end{theorem}
\proof
We proceed to show first the point-wise convergence of $u^{\eps_{k}}$ to $u^{0}$, for some subsequence $\eps_{k} \downarrow 0$, and then, using \Cref{rem:modulustheta}, i.e., constructing such a modulus $\Theta$, we deduce that the convergence is locally uniform. 


From \Cref{thm:spacetightness}, let $\eps_{k}$ be the subsequence such that $m^{\eps_{k}} \to m^{0}$ in $C([0,T];\PP_{1}(\R^{d}))$ as $k \to \infty$. Let $R \geq 0$, let $(t,x,v) \in [0,T] \times \R^{d} \times \overline{B}_{R}$ and let $\gamma^{\eps_{k}}$ be a minimizer for $u^{\eps_{k}}(t,x,v)$. Then, we have that 
\begin{multline*}
u^{\eps_{k}}(t,x,v)
=\ \int_{t}^{T}{\left(\frac{\eps_{k}}{2}|\ddot\gamma^{\eps_{k}}(s)|^{2} + L_{0}(\gamma^{\eps_{k}}(s), \dot\gamma^{\eps_{k}}(s), m^{\eps_{k}}_{s})\right)\ ds} + g(\gamma^{\eps_{k}}(T), m^{\eps_{k}}_{T})
\\
\geq\ \ \int_{t}^{T}{L_{0}(\gamma^{\eps_{k}}(s), \dot\gamma^{\eps_{k}}(s), m^{\eps_{k}}_{s})\ ds}+ g(\gamma^{\eps_{k}}(T), m^{\eps_{k}}_{T})
\\
\geq\  \inf_{\gamma \in \Gamma_{t}(x)}\left\{\int_{t}^{T}{L_{0}(\gamma(s), \dot\gamma(s), m^{0}_{s})\ ds} + g(\gamma(T), m^{0}_{T}) \right\} 
+\  g(\gamma^{\eps_{k}}(T), m^{\eps}_{T}) - g(\gamma^{\eps_{k}}(T), m^{0}_{T})
\\
+\  \int_{t}^{T}{\Big(L_{0}(\gamma^{\eps_{k}}(s), \dot\gamma^{\eps_{k}}(s), m^{\eps_{k}}_{s}) - L_{0}(\gamma^{\eps_{k}}(s), \dot\gamma^{\eps_{k}}(s), m^{0}_{s}) \Big)\ ds}   \geq\  u^{0}(t,x) - o(1)
\end{multline*}
where the last inequality holds by {\bf (M1)} and the convergence of $m^{\eps_{k}}$ in $C([0,T]; \PP_{1}(\R^{d}))$. 

On the other hand, let $R \geq 0$ and take $(t,x,v) \in [0,T] \times \R^{d} \times \overline{B}_{R}$. Let $\gamma^{0} \in \Gamma_{t}(x)$ be a solution of  
\begin{equation*}
	\inf_{\gamma \in \Gamma_{t}(x)} \left\{\int_{t}^{T}{L_{0}(\gamma(s), \dot\gamma(s), m^{0}_{s})\ ds}+g(\gamma(T), m^{0}_{T})\right\}. 
\end{equation*}
Next, we distinguish two cases: first, when $\dot\gamma^{0}(t)=v$ and then when $\dot\gamma^{0}(t)\not=v$. Indeed, if $\dot\gamma^{0}(t)=v$, by the Euler equation and the $C^{2}$-regularity of $L_{0}$ we have that $\gamma \in C^{2}([0,T])$. Hence, we can use $\gamma^{0}$ as a competitor for $u^{\eps_{k}}(t,x,v)$ and we get
\begin{multline}\label{eq:Mfirstest}
u^{\eps_{k}}(t,x,v) \leq  \int_{t}^{T}{\left(\frac{\eps}{2}|\ddot\gamma^{0}(s)|^{2} + L_{0}(\gamma^{0}(s), \dot\gamma^{0}(s), m^{\eps_{k}}_{s})\right)\ ds} + g(\gamma^{0}(T), m^{\eps_{k}}(T))
\\
\leq\ \int_{t}^{T}{\left(\frac{\eps}{2}|\ddot\gamma^{0}(s)|^{2} + L_{0}(\gamma^{0}(s), \dot\gamma^{0}(s), m^{0}_{s})\right)\ ds} + g(\gamma^{0}(T), m^{0}(T))
\\
+\   \int_{t}^{T}{\left(L_{0}(\gamma^{0}(s), \dot\gamma^{0}(s), m^{\eps_{k}}_{s}) - L_{0}(\gamma^{0}(s), \dot\gamma^{0}(s), m^{0}_{s}) \right)\ ds}
\\
+\   g(\gamma^{0}(T), m^{\eps_{k}}(T)) - g(\gamma^{0}(T), m^{0}(T))
\leq\   u^{0}(t,x) + o(1)
\end{multline}
where the last inequality  follows again from the convergence of $m^{\eps_{k}}$ in $C([0,T]; \PP_{1}(\R^{d}))$. If this is not the case, i.e., $\dot\gamma^{0}(t)\not=v$, from \Cref{lem:oscillationvel} we deduce that 
\begin{align*}
u^{\eps_{k}}(t,x,v) =\ u^{\eps_{k}}(t,x,v) - u^{\eps_{k}}(t,x, \dot\gamma^{0}(t)) + u^{\eps_{k}}(t,x,\dot\gamma^{0}(t)) \leq  o(1) + u^{\eps_{k}}(t,x,\dot\gamma^{0}(t)).  	
\end{align*}
Thus, in order to conclude it is enough to estimate $u^{\eps_{k}}(t,x,\dot\gamma^{0}(t))$ as in \eqref{eq:Mfirstest}. Therefore, we obtain
 \begin{equation*}
 	u^{0}(t,x) - o(1) \leq u^{\eps_{k}}(t,x,v) \leq u^{0}(t,x) + o(1)
 \end{equation*}
which implies that $u^{\eps_{k}}$ point-wise converges to $u^{0}$.

Finally, in order to conclude we need to show that the convergence is locally uniform. From \eqref{eq:Mtimegradient}, \eqref{eq:Mspacegradient} and \Cref{lem:oscillationvel} we have that for any $R \geq 0$ and any $(t_{1}, x_{1}, v_{1})$, $(t_{2}, x_{2}, v_{2}) \in [0,T] \times \overline{B}_{R} \times \overline{B}_{R}$ there holds 
\begin{align*}
	& |u^{\eps}(t_{1}, x_{1}, v_{1}) - u^{\eps}(t_{2}, x_{2}, v_{2})|
	\\
	\leq\ & \omega_{R}(|t_{1}-t_{2}|) + C_{1}|x_{1}-x_{2}| + C_{R}\sqrt{\eps}.
\end{align*}
 Therefore, setting 
 \begin{equation*}
 \Theta(\delta_{0}, \varepsilon_{0}, R_{0})=	\omega_{R_{0}}(\delta_{0}) + C_{1}\delta_{0} + C_{R_{0}}\sqrt{\varepsilon_{0}}
 \end{equation*}
by \Cref{rem:modulustheta} we deduce that the convergence is locally uniform and the proof is thus complete. \qed


After proving the convergence of $u^{\eps}$, we go back to the analysis of the flow of measures and in particular we will characterize it in terms of the limit function $u^{0}$. In order to do so, we study the convergence of minimizers for $u^{\eps}$ and appealing to such a result we will show that $m^{0} \in C([0,T]; \PP_{1}(\R^{d}))$ solves a continuity equation with vector field $D_{p}H_{0}(x,D_{x}u^{0})$, in the sense of distribution.

\begin{proposition}\label{prop:convminim}
	Assume {\bf (M1)} -- {\bf (M3)} and {\bf (BC1)}, {\bf (BC2)}. Let $(t, x, v) \in [0,T] \times \R^{2d}$ be such that $u^{0}$ is differentiable at $(t,x)$ and let $\gamma^{\eps}$ be a minimizer for $u^{\eps}(t,x,v)$. Then, $\gamma^{\eps}$ uniformly converges to a curve $\gamma^{0} \in \text{AC}([0,T]; \R^{d})$ and $\gamma^{0}$ is the unique minimizer for $u^{0}(t,x)$ in \eqref{eq:Mlimitvaluefunction}. 
\end{proposition}
\proof
Let us start by proving that $\gamma^{\eps}$ uniformly converges, up to a subsequence. By \Cref{cor:Mapproxbuondedvel} we know that 
\begin{equation*}
\int_{t}^{T}{|\dot\gamma^{\eps}(s)|^{2}\ ds} \leq Q_{1}(1+|v|^{2}).	
\end{equation*}
Thus, for any $s \in [t,T]$, by H\"older inequality we have that 
\begin{align*}
|\gamma^{\eps}(s)| \leq |x| + \sqrt{T}\sqrt{Q_{1}}(1+|v|^{2})^{\frac{1}{2}}. 	
\end{align*}
Therefore, $\gamma^{\eps}$ is bounded in $H^{1}(0,T; \R^{d})$ which implies that by Ascoli-Arzela Theorem there exists a sequence $\{\eps_{k}\}_{k \in \N}$ and a curve $\gamma^{0} \in \text{AC}([0,T]; \R^{d})$ such that $\gamma^{\eps_{k}}$ converges uniformly to  $\gamma^{0}$. 

We show now that such a limit $\gamma^{0}$ is a minimizer for $u^{0}(t,x)$. First, we observe that  
\begin{align*}
&\liminf_{k \to \infty} \left[ \int_{t}^{T}{\left(\frac{\eps_{k}}{2}|\ddot\gamma^{\eps_{k}}(s)|^2 + L_{0}(\gamma^{\eps_{k}}(s), \dot\gamma^{\eps_{k}}(s), m^{\eps}_{s}) \right)\ ds} + g(\gamma^{\eps_{k}}(T), m^{\eps_{k}}_{T})\right] 
\\
\geq\ & \liminf_{k \to \infty} \left[\int_{t}^{T}{L_{0}(\gamma^{\eps_{k}}(s), \dot\gamma^{\eps_{k}}(s), m^{\eps_{k}}_{s})\ ds} + g(\gamma^{\eps_{k}}(T), m^{\eps_{k}}_{T})\right].	
\end{align*}
Then, as observed at the beginning of this proof $\gamma^{\eps}$ is uniformly bounded in $H^{1}(0,T)$. So by lower-semicontinuity of $L$ and \Cref{thm:spacetightness} we deduce that
\begin{align}\label{eq:Mtoprove}
\begin{split}
	& \liminf_{k \to \infty} \left[\int_{t}^{T}{L_{0}(\gamma^{\eps_{k}}(s), \dot\gamma^{\eps_{k}}(s))\ ds} + g(\gamma^{\eps_{k}}(T), m^{\eps_{k}}_{T}) \right]
	\\
	 \geq & \int_{t}^{T}{L_{0}(\gamma^{0}(s), \dot\gamma^{0}(s), m^{0}_{s})\ ds} + g(\gamma^{0}(T), m^{0}_{T}).
	 \end{split}
\end{align}
 Moreover, for any $R \geq 0$ taking $(t,x,v) \in [0,T] \times \R^{d} \times \overline{B}_{R}$, from \Cref{thm:Meps} we obtain 
\begin{equation*}
u^{\eps_{k}}(t,x,v) \leq u^{0}(t,x) + o(1)	
\end{equation*}
 and we recall that 
 \begin{equation*}
 u^{\eps_{k}}(t,x,v) = \int_{t}^{T}{\left(\frac{\eps_{k}}{2}|\ddot\gamma^{\eps_{k}}(s)|^{2} + L_{0}(\gamma^{\eps_{k}}(s), \dot\gamma^{\eps_{k}}(s), m^{\eps_{k}}_{s})\right)\ ds}+g(\gamma^{\eps_{k}}(T), m^{\eps_{k}}_{T}).	
 \end{equation*}
Hence, we get 
\begin{align*}
o(1) + u^{0}(t,x) \geq \int_{t}^{T}{\left(\frac{\eps_{k}}{2}|\ddot\gamma^{\eps_{k}}(s)|^{2} + L_{0}(\gamma^{\eps_{k}}(s), \dot\gamma^{\eps_{k}}(s), m^{\eps_{k}}_{s})\right)\ ds}	+ g(\gamma^{\eps_{k}}(T), m^{\eps_{k}}_{T}).
\end{align*}
Therefore, passing to the limit as $\eps \downarrow 0$ from \eqref{eq:Mtoprove} we obtain 
\begin{equation*}
	u^{0}(t,x) \geq \int_{t}^{T}{L_{0}(\gamma^{0}(s), \dot\gamma^{0}(s), m^{0}_{s})\ ds} + g(\gamma^{0}(T), m^{0}_{T})
	\end{equation*}
	which proves that $\gamma^{0}$ is a minimizer for $u^{0}(t,x)$. Since $u^{0}$  is differentiable at $(t,x) \in \R^{d}$ there exists a unique minimizing trajectory and thus we have that the uniform convergence of $\gamma^{\eps}$ holds for the whole sequence.	 \qed

 \begin{remarks}\em
	Since $u^{0}$ is locally Lipschitz continuous w.r.t. time and Lipschitz continuous w.r.t. space, we have that \Cref{prop:convminim} holds for a.e. $(t,x) \in [0,T] \times \R^{d}$. 
	\end{remarks}


Let $u^{0}$ be as in \eqref{eq:Mlimitvaluefunction} and let $(\gamma^{0}_{t}(\cdot), \dot\gamma^{0}_{t}(\cdot))$ be the flow of Euler-Lagrange equations associated with the minimization problem in \eqref{eq:Mlimitvaluefunction}. Note that, since $u^{0}$ is Lipschitz continuous and $\mu_{0}$ is absolutely continuous w.r.t. the Lebesgue measure we have that on $\supp(\mu_{0})$ the curve $(\gamma^{0}_{t}(\cdot), \dot\gamma^{0}_{t}(\cdot))$ is a minimizer for $u^{0}$. We also recall that the measure $\mu^{\eps}$ is the image of $\mu_{0}$ under the flow \eqref{eq:epsflow}, which is optimal as observed in \Cref{rem:mepsoptimal} for $u^{\eps}(0,x,v)$ for a.e. $(x,v) \in \R^{2d}$, and thus, for any function $\varphi \in C^{\infty}_{c}(\R^{d})$ the measure $m^{\eps}_{t}$ is given by  
\begin{equation}\label{eq:flowconv}
\int_{\R^{d}}{\varphi(x)\ m^{\eps}_{t}(dx)} = \int_{\R^{2d}}{\varphi(\gamma^{\eps}_{(x,v)}(t))\ \mu_{0}(dx,dv)}.	
\end{equation}
We finally recall that by assumption $\mu_{0}$ is absolutely continuous w.r.t. Lebesgue measure.

\begin{corollary}\label{cor:limiteq}
Assume {\bf (M1)} -- {\bf (M3)} and {\bf (BC1)}, {\bf (BC2)}.
Then, we have that
 \begin{equation}\label{eq:zeropush}
m^{0}_{t}=\gamma_{t}^{0}(\cdot) \sharp m_{0}, \quad \forall\ t \in [0,T]. 	
\end{equation}
Moreover, $m^{0} \in C([0,T]; \PP_{1}(\R^{d}))$ solves 
\begin{equation*}
\begin{cases}
\partial_{t} m^{0}_{t} - \ddiv\Big(m^{0}_{t}D_{p}H_{0}(x, D_{x}u^{0}(t,x), m^{0}_{t})\Big)=0, & (t,x) \in [0,T] \times \R^{d}
\\
	m^{0}_{0}=m_{0}, & x \in \R^{d},
	\end{cases}
\end{equation*}
	in the sense of distributions. 
\end{corollary}
\proof
From \Cref{thm:spacetightness} let $\eps_{k} \downarrow 0$ be such that $m^{\eps_{k}} \to m^{0}$ in $C([0,T];\PP_{1}(\R^{d}))$. Then, since $\mu_{0}$ is absolutely continuous w.r.t. Lebesgue measure by \Cref{prop:convminim} we have that 
 \begin{equation*}
	\gamma^{\eps_{k}}_{(x,v)}(t) \to \gamma^{0}_{t}(x), \quad \mu_{0}\text{-a.e.}\, (x,v),\,\, \forall t \in [0,T]. 
\end{equation*}  
Therefore, from \eqref{eq:flowconv}, for $\eps=\eps_{k}$, as $k \to \infty$  we get
\begin{equation*}
	\int_{\R^{d}}{\varphi(x)\ m^{0}_{t}(dx)}= \int_{\R^{d}}{\varphi(\gamma_{t}^{0}(x))\ m_{0}(dx)}, \quad \forall\ t \in [0,T]
\end{equation*}
which proves \eqref{eq:zeropush}.
Moreover, again by \Cref{prop:convminim} we have that $\gamma^{0}_{t}$ is a minimizer for $u^{0}(0,x)$ since it is the limit of $\gamma^{\eps}_{(x,v)}$ which is optimal $u^{\eps}(0,x,v)$ and we are taking $(x,v)$ in a subset of full measure w.r.t. $\mu_{0}$. Therefore, from the optimality of $\gamma^{0}$ we get
\begin{equation*}
\begin{cases}
\dot\gamma^{0}_{t}(x)=D_{p}H_{0}(\gamma_{t}^{0}(x), Du^{0}(t,\gamma_{t}^{0}(x)), m^{0}_{t}), \quad t \in (0,T]
\\
\gamma_{0}^{0}(x)=x.
\end{cases}
\end{equation*}
Hence, for any $\psi \in C^{\infty}_{c}([0,T) \times \R^{d})$ we obtain 
\begin{align*}
& \frac{d}{dt}\int_{\R^{d}}	{\psi(t,x)\ m^{0}_{t}(dx)}=\frac{d}{dt}\int_{\R^{d}}{\psi(t,\gamma_{t}^{0}(x))\ m_{0}(dx)}
\\
= & \int_{\R^{d}}{\big(\partial_{t}\psi(t,\gamma_{t}^{0}(x)) + \langle D_{x}\psi(t, \gamma_{t}^{0}(x)), D_{p}H_{0}(\gamma_{t}^{0}(x), D_{x}\psi(t,\gamma_{t}^{0}(x)), m^{0}_{t})\big)\ m_{0}(dx)}
\\
= & \int_{\R^{d}}{\big(\partial_{t}\psi(t,x) + \langle D_{x}\psi(t, x), D_{p}H_{0}(x, D_{x}\psi(t,x), m^{0}_{t})\big)\ m_{t}^{0}(dx)}
\end{align*}
and integrating, in time, over $[0,T]$ we get the result. \qed
\smallskip

We are now ready to prove the main result. 

\noindent{\it Proof of \Cref{thm:main1}.} Let $\{\eps_{k}\}_{k \in \N}$ be such that $m^{\eps_{k}} \to m^{0}$ in $C([0,T]; \PP_{1}(\R^{d})$ and $u^{\eps_{k}} \to u^{0}$ locally uniformly on $[0,T] \times \R^{2d}$. Then, appealing to \Cref{thm:Meps} and \Cref{cor:limiteq} we deduce that $(u^{0}, m^{0})$ is a solution to the MFG system 
\begin{equation*}
\begin{cases}
	-\partial_{t} u^{0}(t,x) + H_{0}(x, D_{x}u^{0}(t,x), m^{0}_{t})=0, & (t,x) \in [0,T] \times \R^{d} 
	\\
	\partial_{t} m^{0}_{t}- \ddiv\big(m^{0}_{t}D_{p}H_{0}(x, D_{x}u^{0}(t,x), m^{0}_{t})\big)=0, & (t,x) \in [0,T] \times \R^{d}
	\\
	m^{0}_{0}=m_{0}, \,\, u^{0}(T,x)=g(x,m^{0}_{T}) & x \in \R^{d}
\end{cases}	
\end{equation*}
which completes the proof. \qed

		\section{Proof of \Cref{thm:main2}} 
		\label{sec:proof2}

		 We recall that, in this section, we consider the MFG system  
	\begin{align*}
 \begin{cases}
 	-\partial_{t} u^{\eps} +\frac{1}{2\eps}|D_{v}u^{\eps}|^{2} - \langle D_{x}u^{\eps}, v \rangle -\frac{1}{2}|v|^{2} - L_{0}(x, \mu^{\eps}_{t})= 0, & (t,x,v) \in [0,T] \times \R^{2d}
 	\\
 	\partial_{t}\mu^{\eps}_{t} - \langle D_{x}\mu^{\eps}_{t},v \rangle - \frac{1}{\eps}\ddiv_{v}\left(\mu^{\eps}_{t}D_{v}u^{\eps} \right)=0, &  (t,x,v) \in [0,T] \times \R^{2d}
\\
\mu^{\eps}_{0}=\mu_{0}, \quad u^{\eps}(T,x,v)=g(x,\mu^{\eps}_{T}), & (x,v) \in \R^{2d}.
 \end{cases}	
 \end{align*}
 So, the variational problem associated with such a system is given by
 \begin{equation*}
 	u^{\eps}(t,x,v)= \inf_{\gamma \in \Gamma_{t}(x,v)} J^{\eps}_{t,T}(\gamma)
 \end{equation*}
 where
 \begin{align*}
 J^{\eps}_{t,T}(\gamma)=& \int_{t}^{T}{\left(\frac{\eps}{2}|\ddot\gamma(s)|^{2} + \frac{1}{2}|\dot\gamma(s)|^{2} +L_{0}(\gamma(s), \mu^{\eps}_{s}) \right)\ ds} + g(\gamma(T),\mu^{\eps}_{T}), \quad \text{if}\,\, \gamma \in H^{2}(0,T; \R^{d})
 \end{align*}
with $J^{\eps}_{t,T}(\gamma)=+\infty$ if $\gamma \not\in H^{2}(0,T;\R^{d})$. 

From the results on the previous section and the assumptions {\bf (C1)} -- {\bf (C3)} on $L_{0} : \R^{d} \times \PP(\R^{2d}) \to \R$ given above, we deduce that we only need to study the tightness of the flow of measures $\{\mu^{\eps}_{t}\}_{t \in [0,T]}$ w.r.t. the second marginal. This can be done by a finer analysis of the Euler-Lagrange flow. 


\begin{lemma}
Let $(x,v) \in \R^{2d}$ and let $\gamma^{\eps}$ be a solution to the variational problem associated with $u^{\eps}(0,x,v)$. Then, $\gamma^{\eps}$ is a solution of the Euler-Lagrange equation 
\begin{equation*}
-\eps \xi^{(iv)}(t) +\ddot\xi^{\eps}(t) - D_{x}L_{0}(\xi(t), \mu^{\eps}_{t}) = 0
\end{equation*} 
with boundary condition 
\begin{equation*}
\ddot\xi(T)= 0, \quad \text{and} \quad -\eps \xi^{(iii)} (T)  + \dot\xi(T) + D_{x}g(x, \mu^{\eps}_{T}). 
\end{equation*}
\end{lemma}

\begin{proposition}\label{prop:boundacc}
Assume {\bf (C1)} -- {\bf (C3)} and {\bf (A1)}, {\bf (A2)}. Let $(x,v) \in \R^{2d}$ and let $\gamma^{\eps}$ be a solution to the variational problem associated with $u^{\eps}(0,x,v)$. Then, there exists a constant $C \geq 0$ such that for any $\delta \in (0, 1)$ the following holds
\begin{equation*}
\int_{\delta}^{T} |\ddot\gamma^{\eps}(s)|^{2}\ ds \leq C. 
\end{equation*} 
\end{proposition}
\proof

Fix $(x, v) \in \R^{2d}$ and a solution $\gamma^{\eps}$ to the problem associated with $u^{\eps}(0,x,v)$. In the following, for simplicity of notation we drop $\eps$ setting $\gamma^{\eps}=\gamma$ and we drop the notation of the scalar product $\langle \cdot, \cdot \rangle$.


We begin by multiplying the Euler-Lagrange equation 
\begin{equation*}
-\eps \gamma^{(iv)}(t) +\ddot\gamma^{\eps}(t) - D_{x}L_{0}(\gamma(t), \mu^{\eps}_{t}) = 0
\end{equation*} 
by $s^{2}\ddot\gamma(t)$ and integrate by parts. We obtain 
\begin{align*}
& -\eps \gamma^{(iii)} (T)  T \ddot\gamma(T) 
\\
+\ & \int_{0}^{T} \Big(\eps \gamma^{(iii)}(s) \big(s \ddot\gamma(s) + s^{2} \gamma^{(iii)}(s) \big) + s^{2} \ddot\gamma(s)^{2} - D_{x}L_{0}(\gamma(s), \mu^{\eps}_{s}) s^{2}\ddot\gamma(s) \Big)\ ds = 0
\end{align*}
which reduces to 
\begin{equation*}
\int_{0}^{T} \Big(\eps \gamma^{(iii)}(s) \big(s \ddot\gamma(s) + s^{2} \gamma^{(iii)}(s) \big) + s^{2} \ddot\gamma(s)^{2} - D_{x}L_{0}(\gamma(s), \mu^{\eps}_{s}) s^{2}\ddot\gamma(s) \Big)\ ds = 0
\end{equation*}
by using the boundary condition on the Euler-Lagrange equation. From the Young's inequality we get 
\begin{multline*}
 \int_{0}^{T} \Big(\eps \gamma^{(iii)}(s) \big(s \ddot\gamma(s) + s^{2} \gamma^{(iii)}(s) \big) + s^{2} \ddot\gamma(s)^{2}\Big)\ ds 
 \\
\leq\  \int_{0}^{T} \frac{s^{2}}{2}D_{x}L_{0}(\gamma(s), \mu^{\eps}_{s})^{2}\ ds + \int_{0}^{T} \frac{s^{2}}{2}\ddot\gamma(s)^{2}\ ds
\end{multline*}
which yields to 
\begin{equation*}
 \int_{0}^{T} \Big(\eps \gamma^{(iii)}(s) \big(s \ddot\gamma(s) + s^{2} \gamma^{(iii)}(s) \big) + \frac{s^{2}}{2} \ddot\gamma(s)^{2}\Big)\ ds 
\leq\  \int_{0}^{T} \frac{s^{2}}{2}D_{x}L_{0}(\gamma(s), \mu^{\eps}_{s})^{2}\ ds.
\end{equation*}
So, appealing to {\bf (C1)} we have that there exists a constant $C(T) \geq 0$ such that 
\begin{equation*}
\int_{0}^{T} \frac{s^{2}}{2}D_{x}L_{0}(\gamma(s), \mu^{\eps}_{s})^{2}\ ds \leq C(T)
\end{equation*}
and, moreover, by the non-negativity of the term $s^{2}\gamma^{(iii)}(s)^{2}$ we finally get
\begin{equation}\label{eq:b1}
 \int_{0}^{T} \Big(\eps \gamma^{(iii)}(s) s \ddot\gamma(s) + \frac{s^{2}}{2} \ddot\gamma(s)^{2}\Big)\ ds \leq C(T). 
\end{equation}
Integrating, again, by parts the term $ \int_{0}^{T} \eps \gamma^{(iii)}(s) s \ddot\gamma(s)\ ds $ we obtain 
\begin{align*}
\int_{0}^{T} \eps \gamma^{(iii)}(s) s \ddot\gamma(s)\ ds = \eps T \ddot\gamma(T) - \int_{0}^{T} \eps \ddot\gamma(s) \big(\ddot\gamma(s) + s\gamma^{(iii)}(s) \big)\ ds 
\end{align*}
which implies that 
\begin{equation*}
- \int_{0}^{T} \eps \ddot\gamma(s)^{2}\ ds = 2\eps \int_{0}^{T} \gamma^{(iii)}(s) s \ddot\gamma(s)\ ds
\end{equation*}
and so
\begin{equation}\label{eq:b2}
- \int_{0}^{T} \ddot\gamma(s)^{2}\ ds \leq 2\eps \int_{0}^{T} \gamma^{(iii)}(s) s \ddot\gamma(s)\ ds.
\end{equation}
Hence, combining \eqref{eq:b1} and \eqref{eq:b2} we get
\begin{equation*}
\int_{0}^{T} s^{2} \ddot\gamma(s)^{2}\ ds \leq C(T) 
\end{equation*}
which implies the result.
 \qed

\begin{theorem}
Assume {\bf (C1)} -- {\bf (C3)} and {\bf (A1)}, {\bf (A2)}. The sequence $\{\mu^{\eps, \delta}_{t}\}_{\delta > 0}$ is tight and the sequence $\{\mu^{\eps}_{t}\}_{\eps > 0}$ is relatively compact in $C^{0}([0, T]; \PP_{1}(\R^{2d}))$.
\end{theorem}
\proof 
The tightness of the sequence $\{\mu^{\eps, \delta}_{t}\}_{\delta > 0}$ follows from \Cref{cor:Mapproxbuondedvel} and \Cref{prop:boundacc}. 

Next, we show that  $\{\mu^{\eps}_{t}\}_{\eps > 0}$ is relatively compact. Indeed, still from \Cref{cor:Holder} and \Cref{prop:boundacc} we have that for any $\delta > 0$ and any $\delta \leq s \leq t \leq T$
\begin{align*}
d_{1}(\mu_{s}^{\eps}, \mu_{t}^{\eps}) & = \int_{\R^{2d}} \left(|\gamma_{(x,v)}(t) - \gamma_{(x,v)}(s) | + |\dot\gamma_{(x,v)}(t) - \gamma_{(x,v)}(s)| \right)\ \mu_{0}(dx,dv)
\\
\leq\ & (t-s)^{\frac{1}{2}} \left[ \left(\int_{\R^{2d}} |\dot\gamma_{(x,v)}(\tau)|^{2}\ \mu_{0}(dx,dv) \right)^{\frac{1}{2}} + \left(\int_{\R^{2d}} |\ddot\gamma_{(x,v)}(\tau)|^{2}\ \mu_{0}(dx,dv) \right)^{\frac{1}{2}}  \right] 
\\
\leq\ & CQ_{2}(t-s)^{\frac{1}{2}}
\end{align*}
which completes the proof appealing to Prokhorov Theorem and Ascoli-Arzela Theorem. \qed

\begin{remarks}\em
Note that, following the above reasoning one easily deduce that the main result of this section is not uniform w.r.t. $T$. 
\end{remarks}

Now, by using similar techniques of \Cref{thm:Meps} and \Cref{prop:convminim} one can prove the following. 

\begin{proposition}\label{prop:valueconv}
Assume {\bf (C1)} -- {\bf (C3)} and {\bf (A1)}, {\bf (A2)}.  Then, we have that
\begin{itemize}
\item[($i$)] there exists a subsequence $\eps_{k} \downarrow 0$ such that $u^{\eps_{k}}$ locally uniformly converges to $u^{0}$. 
\item[($ii$)] Let $(t,x,v) \in [0,T] \times \R^{2d}$ be such that $u^{0}$ is differentiable at $(t,x)$ and let $\gamma^{\eps}$ be a minimizer of $u^{\eps}(t,x,v)$. Then, $\gamma^{\eps}$ uniformly converges to a curve $\gamma^{0}$ such that $\dot\gamma^{\eps} \to \dot\gamma^{0}$ as $\eps \downarrow 0$ and $\gamma^{0}$ is the unique minimizer for $u^{0}$ in \eqref{eq:Mlimitvaluefunction}.
\end{itemize}
\end{proposition}

\begin{proposition}\label{prop:meaconv}
Assume {\bf (C1)} -- {\bf (C3)} and {\bf (A1)}, {\bf (A2)}. We have that
 \begin{equation*}
\mu^{0}_{t}=(Id_{\cdot}, D_{x}u^{0}(t, \gamma_{t}^{0}(\cdot))) \sharp m^{0}_{t}, \quad \forall\ t \in [0,T]. 	
\end{equation*}
Moreover, $m^{0} \in C([0,T]; \PP_{1}(\R^{d}))$ solves 
\begin{equation*}
\begin{cases}
\partial_{t} m^{0}_{t} - \ddiv\big(m^{0}_{t}D_{x}u^{0}(t,x)\big)=0, & (t,x) \in [0,T] \times \R^{d}
\\
	m^{0}_{0}=m_{0}, & x \in \R^{d},
	\end{cases}
\end{equation*}
	in the sense of distributions. 
\end{proposition}
\proof
From \Cref{prop:valueconv} let $\eps_{k} \downarrow 0$ be such that $\mu^{\eps_{k}} \to \mu^{0}$ in $C([0,T];\PP_{1}(\R^{2d}))$. Then, since $\mu_{0}$ is absolutely continuous w.r.t. Lebesgue measure by \Cref{prop:valueconv} we have that 
 \begin{equation*}
	\gamma^{\eps_{k}}_{(x,v)}(t) \to \gamma^{0}_{t}(x),\,\, \text{and}\,\, \dot\gamma^{\eps_{k}}_{(x,v)}(t) \to \dot\gamma^{0}_{t}(x)  \quad \mu_{0}\text{-a.e.}\, (x,v),\,\, \forall t \in [0,T]. 
\end{equation*}  
Therefore, from \eqref{eq:flowconv}, for $\eps=\eps_{k}$, as $k \to \infty$  we get
\begin{equation*}
	\int_{\R^{d}}{\varphi(x)\ m^{0}_{t}(dx)}= \int_{\R^{d}}{\varphi(\gamma_{t}^{0}(x))\ m_{0}(dx)}, \quad \forall\ t \in [0,T].
\end{equation*}
Moreover, again by \Cref{prop:valueconv} we have that $\gamma^{0}_{t}$ is a minimizer for $u^{0}(0,x)$ since it is the limit of $\gamma^{\eps}_{(x,v)}$ which is optimal $u^{\eps}(0,x,v)$ and we are taking $(x,v)$ in a subset of full measure w.r.t. $\mu_{0}$. Therefore, from the optimality of $\gamma^{0}$ we get
\begin{equation*}
\begin{cases}
\dot\gamma^{0}_{t}(x)= Du^{0}(t,\gamma_{t}^{0}(x)),  \quad t \in (0,T]
\\
\gamma_{0}^{0}(x)=x.
\end{cases}
\end{equation*}
Hence,
 \begin{equation*}
\mu^{0}_{t}=(Id_{\cdot}, D_{x}u^{0}(t, \gamma_{t}^{0}(\cdot))) \sharp m^{0}_{t}, \quad \forall\ t \in [0,T]. 	
\end{equation*}

Then, for any $\psi \in C^{\infty}_{c}([0,T) \times \R^{d})$ we obtain 
\begin{align*}
& \frac{d}{dt}\int_{\R^{d}}	{\psi(t,x)\ m^{0}_{t}(dx)}=\frac{d}{dt}\int_{\R^{d}}{\psi(t,\gamma_{t}^{0}(x))\ m_{0}(dx)}
\\
= & \int_{\R^{d}}{\big(\partial_{t}\psi(t,\gamma_{t}^{0}(x)) + \langle D_{x}\psi(t, \gamma_{t}^{0}(x)), D_{x}\psi(t,\gamma_{t}^{0}(x))\big)\ m_{0}(dx)}
\\
= & \int_{\R^{d}}{\big(\partial_{t}\psi(t,x) + \langle D_{x}\psi(t, x), D_{x}\psi(t,x)\big)\ m_{t}^{0}(dx)}
\end{align*}
and integrating, in time, over $[0,T]$ we get the result. \qed

\noindent{\it Proof of \Cref{thm:main2}.} Let $\{\eps_{k}\}_{k \in \N}$ be such that $\mu^{\eps_{k}} \to \mu^{0}$ in $C([0,T]; \PP_{1}(\R^{2d})$ and $u^{\eps_{k}} \to u^{0}$ locally uniformly on $[0,T] \times \R^{2d}$. Then, appealing to \Cref{prop:valueconv} and \Cref{prop:meaconv} we deduce that $(u^{0}, \mu^{0})$ is a solution to the MFG system 
\begin{align*}
	\begin{cases}
		-\partial_{t} u^{0}(t,x) + \frac{1}{2}|D_{x}u^{0}(t,x)|^{2}-  L_{0}(x, \mu^{0}_{t})=0, & \quad (t,x) \in [0,T] \times \R^{d}
		\\
		\partial_{t}  m^{0}_{t} - \ddiv\big( m^{0}_{t}D_{x}u^{0}(t,x) \big)=0, & \quad (t,x) \in [0,T] \times \R^{d}
		\\
		\mu^{0}_{t} = (\text{Id}_{\cdot}, D_{x}u^{0}(t, \cdot)) \sharp m^{0}_{t}, & \quad t \in [0,T]
		\\
		m^{0}_{0}= m_{0},\,\, u^{0}(T,x)=g(x,\mu^{0}_{T}), & \quad x \in \R^{d},
	\end{cases}
	\end{align*}
which completes the proof. \qed

%
%

\end{document}